\nonstopmode
\input amstex 
\input amsppt.sty   
\hsize 30pc 
\vsize 47pc 
\magnification=\magstep1 
\def\nmb#1#2{\ifx!#1{{\rm(#2)}}\else{#2}\fi} 
\def\cit#1#2{\ifx#1!\cite{#2}\else#2\fi} 
\def\totoc{}             
\def\idx{}               
\def\ign#1{}             
\redefine\o{\circ} 
 
\define\al{\alpha} 
\define\be{\beta} 
\define\ga{\gamma} 
 
\define\ep{\varepsilon}

\define\io{\iota} 
\define\ka{\kappa} 
\define\la{\lambda} 
\define\rh{\rho} 
 
\define\ta{\tau} 
 
\define\ch{\chi} 
\define\ps{\psi} 
 
\define\Ga{\Gamma}

\define\Ph{\Phi} 
\define\Ps{\Psi} 
\define\g{\frak g} 
\define\h{\frak{h}} 
\define\tr{\operatorname{tr}} 
\define\Aut{\operatorname{Aut}} 
\define\End{\operatorname{End}} 
\define\conj{\operatorname{Conj}} 
\define\Ad{\operatorname{Ad}} 
\define\ad{\operatorname{ad}} 
\define\pr{\operatorname{pr}} 
\let\on=\operatorname 
\redefine\i{^{-1}} 
\define\x{\times} 
\def\today{\ifcase\month\or 
 January\or February\or March\or April\or May\or June\or 
 July\or August\or September\or October\or November\or December\fi 
 \space\number\day, \number\year} 
\topmatter 
\title  The generalized Cayley map from an algebraic group to its Lie 
algebra 
\endtitle 
\rightheadtext{The generalized Cayley map} 
\author  Bertram Kostant, Peter W. Michor  \endauthor 
\address B\. Kostant: MIT, Department of Mathematics, 
Cambridge, MA 02139-4307, USA 
\endaddress 
\email kostant\@math.mit.edu \endemail 
\address 
P\. W\. Michor: 
Erwin Schr\"odinger Institute of Mathematical Physics, Boltzmanngasse 
9, A-1090 Wien, Austria, and: 
Institut f\"ur Mathematik, Universit\"at Wien, 
Strudlhofgasse 4, A-1090 Wien, Austria 
\endaddress 
\email Peter.Michor\@esi.ac.at \endemail 
 
\dedicatory For Alexandre Kirillov, on the occasion of his 65-th 
anniversary \enddedicatory 
\thanks 
PWM was supported by `Fonds zur F\"orderung der wissenschaftlichen  
Forschung, Projekt P~14195~MAT.'
\endthanks 
\keywords Cayley mapping for representations \endkeywords 
\subjclass\nofrills{\rm 2000}
 {\it Mathematics Subject Classification}.\usualspace
 Primary 22E50 \endsubjclass 
\abstract 
Each infinitesimally faithful representation of a reductive complex 
connected algebraic group $G$ induces a dominant morphism $\Phi$ from 
the group to its Lie algebra $\g$ by orthogonal projection in the 
endomorphism ring of the representation space. The map $\Phi$ 
identifies the field $Q(G)$ of rational functions on $G$ with an 
algebraic extension of the field $Q(\g)$ of rational functions on 
$\g$. For the spin representation of $\on{Spin}(V)$ the map $\Phi$ 
essentially coincides with the classical Cayley transform. In 
general, properties of $\Phi$ are established and these properties 
are applied to deal with a separation of variables (Richardson) 
problem for reductive algebraic groups: Find $\on{Harm}(G)$ so that 
for the coordinate ring $A(G)$ of $G$ we have $A(G) 
= A(G)^G\otimes \on{Harm}(G)$. As a consequence of a partial solution to 
this problem and a complete solution for $SL(n)$ one has in general 
the equality $[Q(G):Q(\g)] = [Q(G)^G:Q(\g)^G]$ of the degrees of 
extension fields. Among other results, $\Phi$ yields (for the complex 
case) a generalization, involving generic regular orbits, of the 
result of Richardson showing that the Cayley map, when $G$ is 
semisimple, defines an isomorphism from the variety of unipotent 
elements in $G$ to the variety of nilpotent elements in $\g$. In 
addition if $G$ is semisimple the Cayley map establishes a 
diffeomorphism between the real submanifold of hyperbolic elements in 
$G$ and the space of infinitesimal hyperbolic elements in $\g$. Some 
examples are computed in detail. 
\endabstract 
\endtopmatter 
 
\document 
 
 
\head Table of contents \endhead 
\noindent 1. Introduction \leaders \hbox to 
1em{\hss .\hss }\hfill {\eightrm 2}\par  
\noindent 2. The generalized Cayley mapping and its basic properties 
\leaders \hbox to 1em{\hss .\hss }\hfill {\eightrm 4}\par  
\noindent 3. A separation of variables theorem for reductive 
algebraic groups \leaders \hbox to 1em{\hss .\hss }\hfill {\eightrm 10}\par  
\noindent 4. The behavior of the Jordan decomposition under the 
Cayley map \leaders \hbox to 1em{\hss .\hss }\hfill {\eightrm 13}\par  
\noindent 5. The mapping degree of the Cayley mapping 
\leaders \hbox to 1em{\hss .\hss }\hfill {\eightrm 19}\par  
\noindent 6. Examples \leaders \hbox to 
1em{\hss .\hss }\hfill {\eightrm 22}\par  
\noindent 7. Spin representations and Cayley transforms 
\leaders \hbox to 1em{\hss .\hss }\hfill {\eightrm 23}\par  
\noindent References \leaders \hbox to 
1em{\hss .\hss }\hfill {\eightrm 32}\par

\head\totoc\nmb0{1}. Introduction \endhead 
 
Let $G$ be a connected complex reductive algebraic group and let 
$\g=\on{Lie}(G)$. To any rational locally faithful representation 
$\pi:G\to\Aut(V)$ one can associate a dominant morphism 
$$ 
\Ph:G\to \g 
$$ 
which (see example \therosteritem{\nmb|{3}} below) we refer to as a 
generalized Cayley map. If $\pi':\g\to \End(V)$ is the 
differential of $\pi$ the bilinear form 
$(\al,\be)=\tr(\al\be)$ on 
$\End(V)$ defines a projection $pr_{\pi}:\End(V)\to \pi'(\g)$. 
The generalized Cayley map arises from the restriction of $\pr_{\pi}$ 
to $\pi(G)$. This paper establishes a number of striking properties 
of the map $\Ph$. 
 
The map $\Ph$ is conjugation-equivariant and consequently 
{\it $\Ph$ carries conjugacy classes in $G$ to adjoint 
orbits in $\g$} and the corresponding cohomomorphism 
$\Ph^*$ defines a conjugation-equivariant injection 
$$ 
0\longrightarrow A(\g)\longrightarrow A(G) 
$$ 
of affine rings. On the level of quotient fields $\Ph^*$ defines 
$Q(G)$ as a finite algebraic extension of $Q(\g)$. We will write 
$\deg\,\pi$ for 
the degree of this extension. 
 
\subhead Examples \endsubhead 
\roster 
\item"(\nmb:{1})" 
     $G = Gl(n,\Bbb C)$, $\pi$ is the standard representation, 
     $\deg(\pi)=1$. 
\item"(\nmb:{2})" 
     $G = Sl(n,\Bbb C)$, $\pi$ is the standard representation, 
     $\deg(\pi)=n$, see \nmb!{6.2}. 
\item"(\nmb:{3})" 
     $G = \on{Spin}(n,\Bbb C)$, $\pi$ is the Spin representation, 
     $\deg(\pi)=n$ for $n$ even and $\deg(\pi)=n-1$ for $n$ odd, see 
     \nmb!{7.15}. 
\endroster 
 
\proclaim{Theorem} \nmb!{7.15} 
In the above example \therosteritem{\nmb|{3}} the map $\Ph$ 
is, on a Zariski open subset and up to scalar multiplication, given 
by the Cayley transform. 
\endproclaim 
 
Let $A(G)^G$ and $A(\g)^G$ be the subalgebras of $G$-invariants in 
$A(G)$ and $A(\g)$ respectively. The space 
$\on{Harm}(g)\subset A(\g)$ of harmonic polynomials on $\g$ was 
defined in \cit!{9} and the decomposition 
$$ 
A(\g) = A(\g)^G\otimes \on{Harm}(\g) \tag{\nmb:{4}} 
$$ 
was proved in \cit!{9}. As a generalization of \thetag{\nmb|{4}}, 
using Quillen's proof of a conjecture of Serre, Richardson 
in \cit!{16} proved that a $G$-stable subspace $H$ of $A(G)$ exists 
such that 
$$ 
A(G)=A(G)^G\otimes H \tag{\nmb:{5}} 
$$ 
holds. He also raised the question as to 
whether one can give an explicit construction of such a subspace $H$ 
along the lines of \thetag{\nmb|{4}}.  Although this problem is not 
solved in the present paper we do solve a weakened version of the 
problem. Let $\pi$ be given and let $\on{Sing}_{\pi}G$ be the 
subvariety of $G$ where the differential $d\Ph$ is not invertible. 
Then $\on{Sing}_{\pi}G$ is a hypersurface given as the zero set of 
the Jacobian $\Psi\in A(G)^G$ of the mapping $\Ph$, see \nmb!{2.1.2}. 
Let $\on{Harm}(G)=\Ph^*(\on{Harm}(\g))$. We localize with respect to 
$\Psi$ and prove the following 
 
\proclaim{Theorem} \nmb!{3.2} 
One has 
$$ 
A(G)_{\Psi} = A(G)^G_{\Psi}\otimes \on{Harm}(G). \tag{\nmb:{6}} 
$$ 
The statement without localization,
$$ 
A(G) = A(G)^G\otimes \on{Harm}(G), \tag{\nmb:{7}} 
$$ 
holds if and only if $\Ph$ maps regular orbits to regular orbits.
This is the case in example \therosteritem{\nmb|{2}}.
\endproclaim 
 
Moreover \thetag{\nmb|{6}} readily induces the following 
 
\proclaim{Corollary} \nmb!{3.3} 
For the $G$-equivariant 
extension of the rational function fields $\Ph^*:Q(\g)\to Q(G)$ the 
degrees satisfy 
$$ 
[Q(G):Q(\g)] = [Q(G)^G:Q(\g)^G]. 
$$ 
\endproclaim

Let $U\subset G$ be the unipotent variety and $N\subset \g$ the nilcone. 
If $G$ is semisimple and $\pi$ is suitable then Richardson and 
Bardsley in \cit!{1} used $\Ph$, even in the finite characteristic 
case (for good primes), to establish that 
$$ 
\Ph:U\to N \tag{\nmb:{8}} 
$$ 
is an isomorphism of algebraic varieties. We consider here the complex 
case and generalize this for reductive algebraic groups to 

\proclaim{Theorem} \nmb!{4.5}
Let $a\in G$ be regular and assume that $d\Ph(a)$ is invertible.
Then $\Phi$ restricts to an isomorphism 
of affine varieties
$$ 
\Phi:\overline {\conj_G(a)}\to \overline {\Ad_G(\Ph(a))}.\tag{\nmb:{9}} 
$$ 
\endproclaim 
Note that \therosteritem{\nmb|{8}} is 
certainly not an isomorphism in Example \therosteritem{\nmb|{1}}. 

Any 
$a\in G$ has a  
(multiplicative) Jordan decomposition $a=a_sa_u$ where $a_s$ and 
$a_u$ are respectively the semisimple and unipotent components of 
$a$. Analogously any $X\in\g$ has an (additive) Jordan decomposition 
$X= X_s + X_n$ where $X_s$ and $X_n$ are the semisimple and nilpotent 
components. In contrast to \thetag{\nmb|{8}}  the map $\Ph$ always 
carries semisimple elements to semisimple elements. In fact one has 
 
\proclaim{Theorem} \nmb!{4.11} 
For any $a\in G$ one has 
$$ 
\Ph(a_s) = \Ph(a)_s. \tag{\nmb:{10}} 
$$ 
\endproclaim 
 
An element $a\in G$ is called elliptic (resp. hyperbolic) if $a$ 
is semisimple and the eigenvalues of $\pi(a)$ are of norm 1 
(resp. positive) for all $\pi$. Expanding the multiplicative Jordan 
decomposition, every element $a\in G$ has a unique decomposition 
$$ 
a = a_ea_ha_u \tag{\nmb:{11}} 
$$ 
where $a_e$ and $a_h$ are respectively elliptic and hyperbolic and 
all three components commute. We will say that $a$ is of 
\idx{\it positive type} if $a_e$ is the identity. Let 
$G_{\text{pos}}$ be the space of all elements of positive type. 
Analogously an element $X\in \g$ is called elliptic (resp. 
hyperbolic) if $X$ is semisimple and the eigenvalues of $\pi'(X)$ are 
pure imaginary (resp. real) for all $\pi$. Expanding the additive 
Jordan decomposition every element $X\in \g$ has a unique 
decomposition 
$$ 
X = X_e + X_h + X_n \tag{\nmb:{12}} 
$$ 
where $x_e$ and $x_h$ are respectively 
elliptic and hyperbolic and all three components commute. We will say 
that $X$ is of real type if $X_e = 0$. Let $\g_{real}$ be the space 
of all elements of real type in $\g$. Given $\pi$ let 
$G_{\text{nonsing}}$ be the (Zariski open) complement of 
$G_{\text{sing}}$ in $G$. 
 
\proclaim{Theorem} \nmb!{5.5} 
If $a\in G$ is hyperbolic then $\Ph(a)$ is hyperbolic. Furthermore one has 
$G_{\text{pos}}\subset G_{\text{nonsing}}$ and, extending 
\therosteritem{\nmb|{8}}, one has 
$\Ph(G_{\text{pos}})\subset \g_{\text{real}}$ and in fact 
$$ 
\Ph:G_{pos}\to \g_{\text{real}}\tag{\nmb:{13}} 
$$ 
is a diffeomorphism. In particular $\Ph$ defines a bijection of the 
set of all hyperbolic elements in $G$ to the set of all hyperbolic 
elements in $\g$. 
\endproclaim 
 
Obvious questions arise with regard to the restriction of the 
generalized Cayley map to subgroups of $G$. In this connection one 
readily establishes 
 
\proclaim{Theorem} \nmb!{2.6} \nmb!{2.7}
We have 
$$ 
\Ph_{\pi}\vert K = \Ph_{\pi\vert K} \tag{\nmb:{14}}
$$ 
in the following cases:
If $K$ is any reductive subgroup of a reductive $G$ which is the 
connected centralizer of a subset $A\subset G$. Or if $K$ is 
a subgroup corresponding to a simple ideal in the Lie algebra of a 
semisimple group $G$. 
\endproclaim 
 
\subhead Thanks \endsubhead 
We thank N\. Wallach, V.L\. Popov, T.A\. Springer, D\. Luna, D.V\. 
Alekseevsky, J\. Hofbauer, and E\. Vinberg for helpful comments. 
 
\head\totoc\nmb0{2}. The generalized Cayley mapping and its basic 
properties \endhead 
 
\subhead\nmb.{2.1}. The setting\endsubhead 
Let $G$ be a (real or complex) Lie group with 
Lie algebra $\g$, and let $\pi:G\to\Aut(V)$ be a finite dimensional 
representation with $\pi':\g\to \End(V)$ injective. We say that 
\idx{\it $G$ admits a Cayley mapping} if the 
inner product $(A,B)\mapsto \tr(AB)$ on $\End(V)$ 
restricts to a non-degenerate inner 
product $B_\pi$ on $\pi'(\g)$. Thus the orthogonal projection 
$\pr_\pi:\End(V)\to \pi'(\g)$ is well defined. 
We consider the real analytic mapping 
$$\gather 
\Ph=\Ph_\pi:G\to \g\quad\text{ given by }\quad \pi'\o \Ph_\pi 
= \pr_\pi\o\pi:G\to\End(V)\to\pi'(\g),\tag{\nmb:{1}}\\ 
B_\pi(\Ph_\pi(g),X) = \tr(\pi'(\Ph_\pi(g))\pi'(X)) = 
\tr(\pi(g)\pi'(X)), \qquad g\in G, X\in \g, 
\endgather$$ 
which we call the \idx{\it generalized Cayley mapping} 
of the representation $\pi$. The choice of the name is justified by 
the fact that for the Spin representation of $\on{Spin}(n)$ on 
$\Bbb C^n$ it coincides, up to a scalar, 
with the Cayley transform, see \nmb!{7.15}. 
 
By $d\Ph = \pr_2\o T\Ph:TG\to T\g=\g\x\g\to \g$ we denote the 
differential of $\Ph$,  
so that $d\Ph(g)=\pr_2\o T_g\Ph:T_gG\to\g$.  
Let $X_1,\dots,X_n$ be a linear basis of $\g$ and let 
$L_{X_1},\dots,L_{X_n}$ be the corresponding left invariant vector 
fields on $G$. Let 
$$ 
\Ps_\pi(g)=\Ps(g) := \det(d\Ph(g)),\qquad g\in G 
\tag{\nmb:{2}}$$ 
be the \idx{\it Cayley determinant function} of the representation $\pi$, 
where the determinant is computed with respect to the bases 
$L_{X_i}(g)$ of $T_gG$ and $X_i$ of $\g$, respectively. 
Note that $\Ps$ does not depend on the choice of the basis $(X_i)$ of 
$\g$. We get the function $\Ps$ multiplied by the modular function 
if we choose right invariant vector fields. 
 
In the following we shall use the notation $\mu:G\x G\to G$ for the 
multiplication, $\mu(g,h)=gh=\mu_g(h)=\mu^h(g)$ for left and right 
translations, $T(\mu_g)$ and $T(\mu^h)$ for the corresponding tangent 
mappings. 
 
\proclaim{\nmb.{2.2}. Proposition} 
In the following cases the infinitesimally faithful representation 
$\pi:G\to\Aut(V)$ admits a Cayley mapping: 
\roster 
\item"(\nmb:{1})" if $G$ is a reductive complex Lie group 
       and $\pi$ is a holomorphic representation. 
\item"(\nmb:{2})" if $G$ is a real reductive Lie group and $\pi$ is a 
       real or complex representation which maps each element in the 
       connected center to a semisimple transformation (in the 
       complexification of $V$).
       In particular if $G$ is a real compact Lie group. 
\endroster 
\endproclaim 
 
By a real reductive Lie group we mean one where the complexification of 
the Lie algebra is reductive. For abelian Lie groups there are 
representations acting by unipotent matrices only, and we have to 
exclude these. 
 
\demo{Proof} 
\therosteritem{\nmb|{1}} A connected reductive complex Lie group $G$ 
has a compact real form, so the Lie algebra $\g$ of $G$ is the 
complexification of the Lie algebra $\frak k$ of a maximal compact 
subgroup $K\subset G_o$. Since $\pi':\g\to \End(V)$ is complex linear 
it suffices to show that the trace form is non-degenerate on 
$\pi'(\frak k)$. Let us choose a $K$-invariant  
Hermitian product $(\quad,\quad)$ on $V$ by integration. Then $\pi'(X)$ 
is skew Hermitian with respect to this inner product for 
$X\in \frak k$, so 
$0\ge (\pi'(X)v,\pi'(X)v) = -(\pi'(X)^2v,v)$ and $\pi'(X)^2$ is 
negative definite Hermitian, so $\tr(\pi'(X)^2)$ is the sum of the 
negative eigenvalues of $\pi'(X)^2$. Thus $B_\pi$ is non-degenerate 
on $\frak k$. 
 
\therosteritem{\nmb|{2}}
The trace form is non-degenerate on the semisimple part of $\g$, and 
on the center since it is mapped to semisimple endomorphisms.  
For a compact group one can repeat the argument from the proof of 
\therosteritem{\nmb|{1}} with a $G$-invariant positive definite inner 
product on $V$. 
\qed\enddemo 
 
\subhead\nmb.{2.3}. Remark \endsubhead 
Most of the time (when not stated explicitly otherwise) $G$ will 
denote a connected reductive complex algebraic group and $\pi$ will 
be a rational representation; in particular in sections \nmb!{4}, 
\nmb!{3}, and \nmb!{7} below. 
 
Note that the center $\frak c$ of the 
Lie algebra $\g$ then belongs to any Cartan subalgebra, and its Lie 
group (the connected center) to the corresponding Cartan subgroup. 
 
\proclaim{\nmb.{2.4}. Proposition} 
Let $G$ be a complex or real Lie group and let $\Ph$ be the 
generalized Cayley mapping of a representation. 
The Cayley mapping $\Ph$ has the following simple properties: 
\roster 
\item"(\nmb:{1})" $\Ph\o \conj_b = \Ad_b\o\Ph$. 
\item"(\nmb:{2})" For all $g\in G$ we have (where $Z_\g(\g^g)$ 
       denotes the centralizer of $\g^g = \{X\in \g: \Ad_g(X)=X\}$ in $\g$) 
$$\alignat2 
G^g&\subseteq G^{\Ph(g)}, &\qquad \g^a &\subseteq  \g^{\Ph(g)}, \\ 
\Ph(G^g)&\subseteq \g^g, &\qquad 
     \Ph(g) &\in \on{Cent}(\g^g) \subset Z_\g(\g^g). 
\endalignat$$ 
       Moreover, if $G$ is a connected reductive complex algebraic 
       group then we even have $\on{Cent}(\g^g) = Z_\g(\g^g)$. 
\item"(\nmb:{3})" $d\Ph(e):\g\to \g$ is the identity mapping. Thus 
       $d\Ph(g)$ is invertible for $g$ in the non-empty Zariski open 
       dense subset $\{h\in G: \Ps(h)\ne 0\}$ of $G$, and is not 
       invertible on the hypersurface $\{h\in G:\Ps(h)=0\}$. 
\item"(\nmb:{4})" $\Ps$ is invariant under conjugation. 
\item"(\nmb:{5})" Let $H\subset G$ be a Cartan subgroup with Cartan Lie 
     subalgebra $\h\subset \g$. Then $\Ph(H)\subset \h$. 
\item"(\nmb:{6})" Let $\ch_\pi$ be the character of the 
       representation $\pi$, given by $\ch_\pi(g)=\tr(\pi(g))$. Then 
       $d\ch_\pi(g)(T_e(\mu_g)X)= 
       \tr(\pi'(\Ph_\pi(g))\pi'(X))=B_\pi(\Ph_\pi(g),X)$. 
\item"(\nmb:{7})" The differential $d\Ph_\pi(g).T_e(\mu_g).X\in \g$ 
       is given by the implicit equation 
       $\tr(\pi'(d\Ph_\pi(g)\,T_e(\mu_g)\,X)\pi'(Y)) 
          =\tr(\pi(g)\pi'(X)\pi'(Y))$ for $Y\in \g$. 
\item"(\nmb:{8})" If $\Ph(e)=0$ and $a\in G$ is such that 
       $\pi(a)\in \pi'(\g)$ then $d\Ph(a\i)$ is not invertible. 
\endroster 
\endproclaim 
 
\demo{Proof} 
\therosteritem{\nmb|{1}} follows by the invariance of the trace.

\therosteritem{\nmb|{2}} Most of it is obvious.
Let $G$ be a connected reductive complex 
algebraic group. We claim that then $g$ lies in the identity component 
of $G^g$; this is immediate if $g$ is semisimple (see \nmb!{4.1} 
below) since then $g$ lies in a Cartan subgroup. Using the Jordan 
decomposition (see \nmb!{2.6} below) it is true in general. But then 
if $X\in \g$ commutes with $\g^g$ it commutes with $g$, hence 
$X\in \g^g$ and consequently $X\in\on{Cent}(\g^g)$, establishing 
$Z_\g(\g^g)=\on{Cent}(\g^g)$. 
 
\therosteritem{\nmb|{3}} $T_e\pi(G) = \pi'(\g)$. 
 
\therosteritem{\nmb|{4}} follows from \therosteritem{\nmb|{1}} and 
the form of the determinant \nmb!{2.1.2}.
 
\therosteritem{\nmb|{5}} 
Let $a\in H$ be regular in $G$ so that $G^a=H$ and $\g^a=\h$. Then 
use \therosteritem{\nmb|{2}}. 
 
\therosteritem{\nmb|{6}} 
Insert the definitions.
 
\therosteritem{\nmb|{7}} This follows from 
$$\align
\tr(\pi'(d\Ph_\pi(g).&T_e(\mu_g).X)\pi'(Y)) = 
     \tfrac{d}{dt}|_0\tr(\pi'(\Ph_\pi(g\exp(tX)))\pi'(Y)) \\
&= \tfrac{d}{dt}|_0\tr(\pi(g\exp(tX))\pi'(Y))  
= \tfrac{d}{dt}|_0\tr(\pi(g)\pi(\exp(tX))\pi'(Y)) \\
&= \tr(\pi(g)\pi'(X)\pi'(Y)).  
\endalign$$ 
 
\therosteritem{\nmb|{8}} Take $X\in \g$ with $\pi(a)=\pi'(X)$ so that 
$\pi(a\i)\pi'(X) = \on{Id}_V$. Then we have 
$\tr(\pi(a\i)\pi'(X)\pi'(Y))=\tr(\pi'(Y))=0$ for all $Y\in \g$ so 
that $d\Ph(a\i)$ has a non-trivial kernel by 
\therosteritem{\nmb|{7}}. 
\qed\enddemo 
 
\proclaim{\nmb.{2.5}. Proposition} 
Let $G$ be a (real or complex) Lie group and let 
$\pi:G\to\Aut(V)$ be a representation which admits a 
Cayley mapping. 
\roster 
\item"(\nmb:{1})" For $a\in G$ we have 
$$\align 
d\Ph(a)(T_a(\conj_G(a)))&=T_{\Ph(a)}(\Ad_G(\Ph(a)))=[\g,\Ph(a)],\\ 
d\Ph(a)(T_a(G^a))&\subseteq \g^a\subseteq \g^{\Ph(a)},\\ 
T_a(\conj_G(a)) &= \{T_e(\mu^a)X-T_e(\mu_a)X:X\in\g\}\\ 
&= \{T_e(\mu^a)(X-\Ad_aX):X\in\g\},\\ 
T_a(G^a) &= T_e(\mu_a)\g^a = T_e(\mu^a)\g^a.\\ 
\endalign$$ 
      Moreover, if $G$ is a reductive complex algebraic group and $a$ is 
      a semisimple element in $G$ then we have  
$$ 
T_aG = T_a(\conj_G(a))\oplus T_a(G^a),\qquad 
\g= [\g,\Ph(a)] \oplus \g^{\Ph(a)}. 
$$ 
\item"(\nmb:{2})" 
     Suppose that $d\Ph(a):T_aG\to \g$ is invertible for $a\in G$. 
     Then $\g^a=\g^{\Ph(a)}$. In particular, the $G$-orbit 
     $\conj_G(a)$ of $a$ in $G$ has the same dimension as its 
     $\Ph$-image, $\Ad_G(\Ph(a))$; consequently, 
     $\Ph:\conj_G(a)\to \Ad_G(\Ph(a))$ is a covering. If in addition 
     $G$ is connected reductive and $a$ is semisimple, then the 
     centralizer $G_{\Ph(a)}$ is connected, so $\Ph$ is a 
     diffeomorphism between the orbits.
\endroster 
\endproclaim 
 
\demo{Proof} 
\therosteritem{\nmb|{1}} 
Most of it follows easily from \nmb!{2.4}.
Let us now suppose $G$ is a reductive complex algebraic group and 
that $a\in G$ is semisimple \nmb!{4.1} so that $a$ is  
contained in a Cartan subgroup with Cartan Lie subalgebra $\h$. By 
\nmb!{2.4.5} we get $\Ph(a)\in \h$. We claim that then 
$\g= [\g,\Ph(a)] \oplus \g^{\Ph(a)}$. By dimension it suffices to 
check that $[\g,\Ph(a)]\cap \g^{\Ph(a)}=0$. This follows from the 
root space decomposition: Put $\Ph(a)=H\in \h$, let 
$Y=[X,H]\in \g^H$, let $X=X_0+\sum_{\al\in R} X_\al$ be the root 
space decomposition. Then $Y=[X,H]=0+\sum_{\al\in R}\al(H)X_\al$ and 
$0=[Y,H]=\sum_{\al\in R}\al(H)^2X_\al$ whence either $\al(H)=0$ or 
$X_\al=0$ so that $Y=0$. 
 
Similarly it suffices to see that $\g^a\cap (\on{Id}-\Ad_a)\g = 0$. So 
let $a=\exp(H)$ for some $H\in \h$ and suppose that 
$Y=X-\Ad_aX\in \g^a$. Consider again the root space decomposition 
$X=X_0+\sum_{\al\in R} X_\al$. Then we get 
$Y=X-\Ad_{\exp(H)}X=0+\sum_{\al\in R}(1-e^{\al(H)})X_\al$ and 
$0=Y-\Ad_{\exp(H)}Y=\sum_{\al\in R}(1-e^{\al(H)})^2X_\al$ whence 
either $1-e^{\al(H)}=0$ or $X_\al=0$ so that $Y=0$. 
 
\therosteritem{\nmb|{2}} 
By \nmb!{2.4.2} we have $\g^a\subseteq \g^{\Ph(a)}$, so it suffices 
to prove the converse. Let $X\in \g^{\Ph(a)}$ so that $[X,\Ph(a)]=0$. 
Then 
$$\align 
d\Ph(a)T_e(\mu^a)(X-\Ad_aX) &= \tfrac{d}{dt}|_0 \Ph(\exp(tX)a\exp(-tX))\\ 
&= \tfrac{d}{dt}|_0 \Ad_{\exp(tX)}\Ph(a) = [X,\Ph(a)] = 0\\ 
\endalign$$ 
so that $X-\Ad_aX=0$ since $d\Ph(a)$ is invertible. 
If $G$ is connected reductive and $a$ is semisimple (see \nmb!{4.1}) 
then by \nmb!{2.4.5} $\Ph(a)$ is also semisimple, thus $G^{\Ph(a)}$ 
is connected. 
\qed\enddemo 
 
\proclaim{\nmb.{2.6}. Theorem} 
Let $G$ be a (real or complex) Lie group and let 
$\pi:G\to\Aut(V)$ be a representation which admits a 
Cayley mapping. 
Let $H=(\bigcap_{a\in A}G^a)_o=(G^A)_o\subseteq G$ be a subgroup 
which is the connected centralizer of a subset $A\subseteq G$ and 
suppose that $H$ is itself reductive. 
 
Then the representation $\pi|H:H\to \End(V)$ admits a Cayley mapping and 
we have 
$$ 
\Ph_\pi|H = \Ph_{\pi|H}:H \to \h. 
$$ 
\endproclaim 
 
Note that this result proves again \nmb!{2.4.5} by choosing 
$A$ the Cartan subgroup. If we choose $A=G$ then $H$ is the connected 
center of $G$.  
The assumption of this theorem holds in the following two cases: 
\roster 
\item"(\nmb:{1})" If $H=(G^g)_o$ is the connected centralizer of a 
       semisimple element $g$ in a complex algebraic group $G$ (see 
       \nmb!{4.1} below), since then  
       the connected centralizer $(G^g)_o=Z_G(g)_o$ is  
       again reductive, by Proposition 13.19 in \cit!{2} on p\. 321. 
\item"(\nmb:{2})" If $H=(G^A)_o$ is the connected centralizer of a 
       reductive subgroup $A\subseteq G$ of a complex reductive Lie 
       group $G$, since it is well known that 
       then $H$ is itself reductive. 
\endroster

\demo{Proof} 
By \nmb!{2.2.1} both Cayley mappings $\Ph_\pi$ and 
$\Ph_{\pi|H}$ exist. 
By \nmb!{2.1.2} we have $\Ph_\pi(G^a)\subseteq \g^a$ for all 
$a\in A$, thus also 
$$ 
\Ph_\pi(H)\subseteq \bigcap_{a\in A}\Ph_\pi(G^a)\subseteq 
 \bigcap_{a\in A}\g^a=\h. 
$$ 
Moreover, for $h\in H$ and $X\in \h$ we have 
$$ 
\tr(\pi'(\Ph_\pi(h)\pi'(X))=\tr(\pi(h)\pi'(X)) 
     =\tr(\pi'(\Ph_{\pi|(G^g)_o}(h)\pi'(X))). 
$$ 
Thus $\Ph_\pi(h)=\Ph_{\pi|H}(h)$. 
\qed\enddemo 
 
\proclaim{\nmb.{2.7}. Theorem} 
Let $G$ be a semisimple real or complex Lie group, 
let $\pi:G\to\Aut(V)$ be an infinitesimally effective representation. 
Let $\g=\g_1\oplus\dots\oplus\g_k$ be the decomposition 
into the simple ideals $\g_i$. 
Let $G_1,\dots,G_k$ be the corresponding connected subgroups of $G$. 
Then 
$$ 
\Ph_\pi|G_i = \Ph_{\pi|G_i}\quad\text{ for }i=1,\dots,k. 
$$ 
\endproclaim 
 
This result cannot be extended to reductive groups as the standard 
representation of $GL_n$ shows where $\Ph$ is the embedding 
$GL_n\to \frak g\frak l_n$.  
 
\demo{Proof} 
Let $g_1\in G_1$ and suppose that 
$\Ph_\pi(g_1)=\sum_{i=1}^k\Ph_\pi(g_1)_i\in \bigoplus_{i=1}^k\g_i$ 
with $\Ph_\pi(g_1)_i\ne 0\in \g_i$ for some $i\ne 1$. Then 
for each $g_i\in G_i$ we have 
$$\align 
\Ph_\pi(g_1) &= \Ph_\pi(g_i\,g_1\,g_i^{-1}) = \Ad_{g_i}\Ph_\pi(g_1) 
\endalign$$ 
so that $\Ph_\pi(g_1)$ contains a nontrivial $G_i$-orbit which is absurd. 
This shows that $\Ph_\pi(G_j)\subseteq\g_j$ for $j=1,\dots,k$. 
Moreover for $g_j\in G_j$ and $X_j\in \g_j$ we have 
$$ 
\tr(\pi'(\Ph_\pi(g_j))\pi'(X_j)) = \tr(\pi(g_j)\pi'(X_j)) 
= \tr((\pi|G_j)'(\Ph_{\pi|G_j}(g_j)(\pi|G_j)'(X_j)) 
$$ 
which implies the result. 
\qed\enddemo 
 
\proclaim{\nmb.{2.8}. Theorem} 
Let $G$ be a simple real or complex Lie group, 
and let $\pi_i:G\to \End(V_i)$ be nontrivial  
representations for $i=1,2$. The inner product $B_{\pi_i}$ 
on $\g$ from \nmb!{2.1} is a multiple of the Cartan Killing form $B$, 
we write $B_{\pi_i}=j_{\pi_i}\,B$. Then we have. 
\roster 
\item"(\nmb:{1})" 
       For the direct sum representation  
       $\pi_1\oplus \pi_2:G\to \End(V_1\oplus V_2)$ we have 
$$ 
\Ph_{\pi_1\oplus \pi_2}(g) = 
     \frac{j_{\pi_1}}{j_{\pi_1\oplus \pi_2}}\,\Ph_{\pi_1}(g) 
     + \frac{j_{\pi_2}}{j_{\pi_1\oplus \pi_2}}\,\Ph_{\pi_2}(g)  \in \g. 
$$ 
\item"(\nmb:{2})" 
       For the tensor product representation 
       $\pi_1\otimes \pi_2:G\to \End(V_1\otimes V_2)$ we have 
$$ 
\Ph_{\pi_1\otimes \pi_2}(g) = 
   \frac{j_{\pi_1}\ch_{\pi_2}(g)}{j_{\pi_1\otimes \pi_2}}\,\Ph_{\pi_1}(g) 
   + \frac{\ch_{\pi_1}(g)j_{\pi_2}}{j_{\pi_1\otimes \pi_2}}\,\Ph_{\pi_2}(g) 
   \in \g. 
$$ 
\item"(\nmb:{3})" 
        For the $n$-fold tensor product representation 
        $\otimes^n\pi:G\to \End(\otimes^n V)$ we have 
$$ 
\Ph_{\otimes^n\pi}(g) = 
     \left(\frac{\ch_\pi(g)}{\dim(V)}\right)^{n-1}\Ph_\pi(g). 
$$ 
\item"(\nmb:{4})" 
        For the contragredient representation 
        $\pi^\top:G\to \End(V^*)$ given by 
        $\pi^\top(g)=\pi(g\i)^\top$ we have 
$$ 
\Ph_{\pi^\top}(g) = -\Ph_\pi(g\i). 
$$ 
\endroster 
\endproclaim

For a complex simple Lie group the number $j_{\pi_i}$ is a multiple of the 
\idx{\it Dynkin index} of the representation $\pi_i$. It is 
non-negative and satisfies  
$$\align 
j_{\pi_1\oplus\pi_2}&=j_{\pi_1}+ j_{\pi_2},\\ 
j_{\pi_1\otimes\pi_2}&=\dim(V_2)j_{\pi_1}+ \dim(V_1)j_{\pi_2}, 
     \tag{\nmb:{5}}\\ 
j_{\pi}&=\frac{\dim(V)}{\dim(\g)}B(\la_{\pi},\la_{\pi}+\rh), 
\endalign$$ 
for an irreducible representation $\pi$ with highest weight 
$\la_{\pi}$, where $\rh$ is half the sum of all positive roots. This 
is due to Dynkin \cit!{4} and can be found in \cit!{5},~p.100. 
Using this, equation \therosteritem{\nmb|{2}} becomes 
$$\align 
\Ph_{\pi_1\otimes \pi_2}(g) 
&= \frac{B(\la_{\pi_1},\la_{\pi_1}+\rh)} 
     {B(\la_{\pi_1},\la_{\pi_1}+\rh)+B(\la_{\pi_2},\la_{\pi_2}+\rh)} 
     \,\frac{\ch_{\pi_2}(g)}{\dim(V_2)}\,\Ph_{\pi_1}(g) \tag{\nmb|{2}'}\\ 
&\qquad+ \frac{B(\la_{\pi_2},\la_{\pi_2}+\rh)} 
     {B(\la_{\pi_2},\la_{\pi_2}+\rh)+B(\la_{\pi_1},\la_{\pi_1}+\rh)} 
     \,\frac{\ch_{\pi_1}(g)}{\dim(V_1)}\,\Ph_{\pi_2}(g) \\ 
\endalign$$ 
 
\demo{Proof} 
\therosteritem{\nmb|{1}}, \therosteritem{\nmb|{2}} and 
\therosteritem{\nmb|{4}} are easy computations. 

\therosteritem{\nmb|{3}} By induction and \thetag{\nmb|{5}} we check that 
$j_{\otimes^n\pi}=n\dim(V)^{n-1}j_\pi$ which via \therosteritem{\nmb|{2}} 
leads quickly to the result. 
\qed\enddemo 
 
\proclaim{\nmb.{2.9}. Proposition} 
For the Cayley mapping of a rational representation of a connected 
reductive complex algebraic group $G$,  
the pullback mapping $\Ph^*:A(\g)=S^*(\g^*) \to A(G)$ between the 
algebras of regular functions is injective, equivariant, and maps the 
subalgebras of invariant regular functions to each other, 
$\Ph^*:A(\g)^G\to A(G)^G$. Consequently, $\Ph:G\to \g$ is a dominant 
algebraic morphism. 
By the algebraic Peter-Weyl 
theorem we have $A(G)= \bigoplus_{\la\in D} A_\la$ where $D$ is 
the set of all dominant integral highest weights, and where
$$
A_\la=\{f\in A(G): f(g)=\tr(\pi_\la(g)B) \text{ for some 
}B\in \End(V_\la) \}
$$
For an irreducible representation $\pi$ we 
thus have $\Ph^*(\g^*)\subset A_\la$ where $\la$ is the highest 
weight of $\pi$. 
\endproclaim 

\demo{Proof} We prove the last statement. Let $\la$ be the highest 
weight of $\pi$. 
$A_\la$ is a vector space of dimension $(\dim(V_\la))^2$ of functions 
on $G$. 
Take $X\in\g$ and consider the function $f\in A_\la$ given by 
$f(g)=\tr(\pi(g)\pi'(X))$. 
Then $f=\Ph^*(B_\pi(\quad,X))$ 
where $B_\pi(\quad,X)\in\g^*$. 
\qed\enddemo

\head\totoc\nmb0{3}. A separation of variables theorem for reductive 
algebraic groups \endhead 
 
In this section $G$ is a connected reductive complex algebraic group 
and $\Ph$ is the Cayley mapping of a rational 
representation. 
 
\subhead\nmb.{3.1}. Rational vector fields on $G$ \endsubhead 
Let $x_i:\g\to \Bbb C$ be the coordinate functions for the basis 
$X_1,\dots,X_n$ of $\g$, and consider the constant vector fields 
$\partial_{x_i}=\frac{\partial}{\partial x_i}$. 
 
\proclaim{Theorem} 
Let $G$ be a connected reductive complex algebraic group 
and let $\Ph$ be the Cayley mapping of a rational 
representation. 
 
Then the pull back vector fields $Y_i:=\Ph^*\partial_{x_i}$ for 
$i=1,\dots,n$ are commuting vector fields with rational coefficients 
on $G$. The fields $\Ps.Y_i$ are regular vector fields (with 
algebraic coefficients). The fields $Y_i$ induce an $G$-equivariant 
injective algebra homomorphism from the symmetric algebra $S^*\g$ 
into the algebra of differential operators on $G$ with rational 
coefficients with polar divisors contained in the hypersurface 
$\{h\in G:\Ps(h)=0\}$. 
\endproclaim 
 
\demo{Proof} 
For a matrix $A$ we denote by $C(A)$ the classical adjoint or the 
matrix of algebraic complements which 
satisfies $A.C(A)=C(A).A = \det(A).\on{Id}$. Applying this to 
$d\Ph(g)$ we can write the pull back fields as 
$$\align 
(\Ph^*\partial_{x_i})(g) &= d\Ph(g)\i.\partial_{x_i}|_{\Ph(g)} 
     = \frac1{\det(d\Ph(g))}.C(d\Ph(g)).\partial_{x_i}\\ 
&=\frac1{\Ps(g)}.\sum_{j=1}^n C(d\Ph(g))_{ij} L_{X_j}(g) 
\endalign$$ 
These are well defined algebraic vector fields on the Zariski open 
set $\{\Ps\ne 0\}$ and are $\Ph$-related to the constant fields 
$\partial_{x_i}$ on $\g$, thus they commute. 
\qed\enddemo 
 
\subhead\nmb.{3.2}. Invariants and harmonic functions \endsubhead 
For the algebra of regular functions we have 
$A(\g)=A(\g)^G\otimes \on{Harm}(\g)$ by \cit!{9}, theorem 
0.2., where the space $\on{Harm}(\g)$ is by definition the space of 
all regular functions which are killed by all invariant differential 
operators with constant coefficients. 
We define $\on{Harm}_\pi(G):= \Ph_\pi^*(\on{Harm}(\g))$. 
It is a $G$-module. Let us 
denote by $A(G)_{\Ps}$ the localization at $\Ps$. 
 
\proclaim{Theorem} 
Let $G$ be a connected reductive complex algebraic group 
and let $\Ph$ be the Cayley mapping of a rational 
representation. Then 
$$ 
A(G)_{\Ps} = A(G)^G_{\Ps}\otimes \on{Harm}_\pi(G). 
$$ 
Moreover, we have 
$$ 
A(G) = A(G)^G\otimes\on{Harm}_\pi(G)
$$
if and only if $\Ph:G\to \g$ maps regular orbits in $G$ to regular 
orbits in $\g$.
\endproclaim 
 
Note that for the standard representation $\pi$ of $SL_n(\Bbb C)$ the 
Cayley mapping $\Ph_\pi$ carries regular orbits to regular orbits, 
see \nmb!{6.2} below. In general this is wrong. If there exists a 
closed orbit $O=\conj_{G}(a)$ in $G$ for $a$ in a Cartan subgroup $H$, 
such that 
$\dim(\Ph(O))<\dim(O)$ then there exists already a regular orbit with 
that property,  $\Ad(G)(a\,u)\subset \g$ where $u$ is a 
principal unipotent element in $(G^a)_o$, see \nmb!{4.8.3}. 
The orbits through 
$$ 
\on{diag}(1,1,1,1,-1,-1,-1,-1)\in SO(8)\cap \frak s\frak o(8) 
$$ 
have this property.  
 
\demo{Proof} 
By \cit!{15}, theorem A, we have 
$$ 
A(G)= A(G)^G\otimes H 
\tag{\nmb:{1}}$$ 
for some $G$-submodule $H$ of $A(G)$. Then $H$ is by restriction 
isomorphic to the affine ring $A(\conj_G(a))$ for any regular orbit 
$\conj_G(a)$ in $G$. 
 
Let $\la$ be any highest weight of $G$ such that the corresponding 
irreducible $G$-module $V_\la$ has a non-zero weight space $V_\la^0$ 
with weight 0 of dimension $d(\la)$. Let $H_\la$ be the primary 
component of $H$ of type $V_\la$ so that one has the direct sum 
$$H=\bigoplus_\la H_\la.$$ 
By \cit!{9}, p~348, proposition~8, the multiplicity of $V_\la$ in 
$H_\la$ is $d(\la)$, since $H$ restricts bijectively to any regular 
orbit in $G$; see also \cit!{15}, theorem A. 
Thus we may write as a direct sum 
$$ 
H_\la = \bigoplus_{j=1}^{d(\la)} H_\la^j 
\tag{\nmb:{2}}$$ 
where $H_\la^j$ is irreducible and hence equivalent to $V_\la$. 
 
Now let $\on{Harm}(G)_\la$ be the primary component of the 
$G$-module $\on{Harm}(G)=\Ph^*(\on{Harm}(\g))\subset A(G)$ of 
type $V_\la$. By \cit!{9} the multiplicity of $V_\la$ in 
$\on{Harm}(G)$ is again $d(\la)$ so that we have 
$$ 
\on{Harm}(G)_\la = \bigoplus_{j=1}^{d(\la)} \on{Harm}(G)_\la^j 
\tag{\nmb:{3}}$$ 
where each $\on{Harm}(G)_\la^j$ is an irreducible $G$-submodule of 
$A(G)= A(G)^G\otimes H$ of highest weight $\la$ and hence 
equivalent to $V_\la$. Thus there exists a $(d(\la)\x d(\la))$-matrix 
$S=(s_{ij})$ with entries in $A(G)^G$ so that 
$$ 
\on{Harm}(G)_\la^i = \sum_{j=1}^{d(\la)}s_{ij} H_\la^j. 
\tag{\nmb:{4}}$$ 
One would like to replace the abstract $G$-module $H_\la$ by the 
equivalent and explicit $G$-module $\on{Harm}(G)_\la$. For that one 
needs that the matrix $S$ is invertible in $A(G)^G$. The 
determinant $\det(S)$ is non-zero in $A(G)^G$ by the 
independence in \thetag{\nmb|{2}} and \thetag{\nmb|{3}}. 
One would need that $\det(S)$ 
is a constant. Let $Z=\on{Zero(\det(S))}$ be the zero set of 
$\det(S)$. By Steinberg \cit!{18}, theorem~1.3, the set of irregular 
elements in $G$ is of codimension 3. Since $\det(S)\in A(G)^G$ 
its zero set $Z$ contains full orbits and is of codimension 1, so it 
is the union of the closures of all regular orbits  $\Cal O$ in $Z$. 
But $\Ph(\Cal O)$ cannot be a regular adjoint orbit in $\g$ since 
$\on{Harm}(\g)$ restrict faithfully to $A(\Ad(G).X)$ for every 
regular orbit in $\g$. 
Consequently, if one knew that $\Ph$ carried 
each regular orbit to a regular orbit then $Z$ was empty and hence 
$\det(S)$ a non-zero constant which proves part of the second assertion.

In the general case, if we localize $A(G)$ and $A(G)^G$ at 
the function $\Ps$ from \nmb!{2.4}, then we restrict $\Ph$ to the 
Zariski open affine subvariety $G_{\Ps\ne 0}$ where $\Ps$ does not 
vanish. But there $\Ph$ is locally biholomorphic and thus carries 
regular orbits to regular orbits in $\g$, see \nmb!{2.5.2}. 
Hence $\det(S)$ does not 
vanish on $G_{\Ps\ne 0}$ and is thus invertible in $A(G)^G_\ps$. 
Consequently we have $A(G)_\Ps=A(G)^G_\Ps.\on{Harm}(G)$. 
 
It remains to prove that 
$A(G)_\Ps=A(G)^G_\Ps.\on{Harm}(G)$ implies 
$A(G)_\Ps=A(G)^G_\Ps\otimes\on{Harm}(G)$. 
Assume for contradiction that this is false. 
Then there exist linearly independent $b_i\in \on{Harm}(G)$ and 
linearly independent $a_i\in A(G)^G_\Ps$ for $1\le i\le k$ such  
that $\sum_i a_ib_i=0$. Since $G_{\Ps\ne0}$ is Zariski open in 
$G$ there exists a regular orbit $\Cal O$ in $G_{\Ps\ne0}$ such that 
$a_i(g)\ne 0$ for all $i$ and any $g\in \Cal O$. On the other hand 
$\sum_i a_ib_i|\Cal O=\sum_i a_i(g)b_i|\Cal O$ vanishes on $\Cal O$. 
But $\Ph(\Cal O)$ is a regular orbit in $\g$. Let $b_i=\Ph^*(c_i)$ 
for $c_i\in \on{Harm}(\g)$; then the $c_i$ are linearly independent. 
Thus the vanishing of $\sum_i a_i(g) c_i$ on the regular orbit 
$\Ph(\Cal O)$ contradicts the fact that $\on{Harm}(\g)$ restricts 
faithfully onto each regular orbit. 

Finally, if we had $A(G) = \on{Harm}_\pi(G))\otimes A(G)^G$ then 
$\on{Harm}_\pi(G))=\Ph^*\on{Harm}(\g))$ would restrict faithfully to 
each regular orbit in $G$. Since the same is true for 
$\on{Harm}(\g))$ the mapping $\Ph$ had to map regular orbits to 
regular orbits. 
\qed\enddemo 
 
\proclaim{\nmb.{3.3}. Corollary} 
Let $G$ be a connected reductive complex algebraic group 
and let $\Ph$ be the Cayley mapping of a rational 
representation with $\Ph(e)=0\in \g$. 
 
Then for the $G$-equivariant 
extension of the rational function fields $\Ph^*:Q(\g)\to Q(G)$ the 
degrees satisfy 
$$ 
[Q(G):Q(\g)] = [Q(G)^G:Q(\g)^G]. 
$$ 
\endproclaim 
 
See \nmb!{6.2} below for the explicit extension in the case of the standard 
representation of $SL_n(\Bbb C)$. 
 
\demo{Proof} 
Note that $Q(\g)^G$ is the quotient field of $A(\g)^G$, and $Q(G)^G$ 
is the quotient field of $A(G)^G$. We have 
$A(\g)=A(\g)^G\otimes \on{Harm}(\g)$ by \cit!{9}, theorem 0.2., and 
$A(G)_{\Ps} = A(G)^G_{\Ps}\otimes \on{Harm}_\pi(G)$ by \nmb!{3.2}, 
where  $\on{Harm}_\pi(G):= \Ph^*(\on{Harm}(\g))$ is isomorphic to 
$\on{Harm}(\g)$. 
 
Let $k=[Q(G),Q(\g)]$. Then any $q\in Q(G)$ satisfies a unique 
monic polynomial of degree $\le k$ with coefficients in $Q(g)$. 
Choose $q\in Q(G)^G$. Then the coefficients must be in $Q(g)^G$ since 
otherwise by conjugating by an element in $G$ we would obtain a new 
minimal polynomial which contradicts uniqueness. Thus 
$[Q(G)^G,Q(\g)^G]\le k$. 
 
On the other hand if $m=[Q(G)^G,Q(\g)^G]$ then there exists $q$ in 
$Q(G)^G$ which satisfies an equation of degree $m$ over $Q(\g)^G$ and 
$Q(G)^G=(Q(\g)^G)[q]$. But $A(G)^G_{\Psi}$ is contained in $Q(G)^G$. 
Thus $A(G)_{\Psi}$ is contained in $Q(\g)[q]$ by \nmb!{3.2}. Hence 
the quotient field $Q(G)$ of $A(G)_{\Psi}$ is contained in 
$Q(\g)[q]$. Thus $m=k$. 
\qed\enddemo 
 
\head\totoc\nmb0{4}. The behavior of the Jordan decomposition under 
the Cayley map 
\endhead 
 
In this section $G$ is a connected reductive complex algebraic group 
and $\Ph$ is the generalized Cayley mapping of a rational 
representation. 
 
\subhead\nmb.{4.1}. The Jordan decomposition \endsubhead 
For references about the Jordan decomposition (additive in the Lie 
algebra and multiplicative in the algebraic group) see \cit!{2} or 
chapter IV of \cit!{6}. In our case given $a\in G$ we write $a_s$ for 
the semisimple part of $a$ and $a_u$ for the unipotent part of $a$ 
so that $a = a_s\,a_u = a_u\,a_s$.  Recall that semisimple means that 
$a_s$ is $G$-conjugate to an element in the Cartan subgroup $H$ and 
$a_u$ unipotent means 
that $a_u$ is conjugate to an element in the unipotent variety $U$. 
We shall use the decomposition $H=T\,H_{\Bbb R}$ where $T$ is a 
maximal torus.  
An element $a\in G$ is called elliptic if $a$ 
is semisimple and the eigenvalues of $\pi(a)$ are of norm 1 
for all $\pi$; equivalently, $a$ is conjugated to an element in $T$. 
Likewise, an element $a\in G$ is called hyperbolic if $a$ 
is semisimple and the eigenvalues of $\pi(a)$ are real for all $\pi$; 
or equivalently, if $a$ is conjugated to to an element in $H_\Bbb R$. 
Expanding the multiplicative Jordan 
decomposition every element $a\in G$ has a unique decomposition 
$$ 
a = a_ea_ha_u 
$$ 
where $a_e$ and $a_h$ are respectively elliptic and hyperbolic and 
all three components commute. 
We say that $a$ is \idx{\it of positive type} if $a_e=e\in G$. 
 
Analogously for $X\in \g$ there is the unique (additive) Jordan 
decomposition $X = X_s + X_n$ where $[X_s,X_n] = 0$ and where $X_s$ 
is semisimple (conjugate to an element in $\h$) and $X_n$ is 
nilpotent (conjugate to an element in the nilcone $N$). 
Expanding the additive 
Jordan decomposition every element $X\in \g$ has a unique 
decomposition 
$$ 
X = X_e + X_h + X_n 
$$ 
where $X_e$ and $X_h$ are respectively 
elliptic and hyperbolic and all three components commute. We will say 
that $X$ is \idx{\it of real type} 
if $X_e = 0$. Let $\g_{real}$ be the space 
of all elements of real type in $\g$. 
For information on hyperbolic elements in algebraic groups see 
\cit!{11}, especially Section 2 on p\.~418. 
 
\subhead\nmb.{4.2}. $\Ad$-complete modules \endsubhead 
Let $D\subset h^*$ denote the set of dominant 
integral weights for $G$ (relative to some 
fixed Borel subgroup) and for each 
$\la\in D$ let $\pi_{\la}:G\to 
Aut\,V_{\la}$ be a fixed irreducible 
representation with highest weight 
$\la$. 
 
A completely reducible $G$-module $M$ will 
be said to $\Ad$-complete if one has an 
equivalence 
$$ 
M\cong \oplus_{\la\in D}\dim(V_{\la}^H)\,V_{\la}, 
$$ 
i.e.\,  each irreducible component occurs with 
multiplicity equal to the dimension of 
its zero weight space. 
 
Let $\on{Reg}(\g)$ (resp\. $\on{Reg}(G)$) be the set of regular 
elements in $\g$ (respectively $G$. We recall the following results. 
 
\proclaim{\nmb.{4.3}. Theorem} \cit!{9} 
For $X\in \on{Reg}(\g)$ one has 
$A(\Ad_G(X))= A(\overline {\Ad_G(X)})$ 
and as a $G$-module $A(\overline{\Ad_G(X)})$ is $\Ad$-complete. 
Furthermore there exists 
a section $\on{Reg}_{\#}(\g)$ of the map 
$\on{Reg}(\g)\to \on{Reg}(\g)/G$ which in addition has the property that  
$$ 
\on{Reg}_{\#}(\g)\to \Bbb C^{\ell},\qquad X\mapsto 
     (I_1(X),...,I_{\ell}(X))\tag{\nmb:{1}} 
$$ 
is an algebraic isomorphism, where the $I_k$ form a basis of $A(\g)^G$. 
Finally, 
$$ 
\g = \bigcup_{X\in \on{Reg}_{\#}(\g)}{\overline {\Ad_G(X)}}\tag{\nmb:{2}} 
$$ 
is a disjoint union. 
\endproclaim 
 
Subsequently Steinberg proved the following group-theoretic analogue: 
 
\proclaim{\nmb.{4.4}. Theorem} \cit!{18} 
Assume $G$ is simply-connected semisimple. 
For $a\in \on{Reg}(G)$ one has 
$$ 
A(\conj_G(a))= A(\overline {\conj_G(a)}) 
$$ 
and as a $G$-module $A(\overline {\conj_G(a)})$ is $\Ad$-complete. 
Furthermore there exists 
a section $\on{Reg}_{\#}(G)$ of the map 
$\on{Reg}(G)\to \on{Reg}(G)/G$ which in addition has the property that  
$$ 
\on{Reg}_{\#}(G)\to \Bbb C^{\ell},\qquad a\mapsto 
(\chi_1(a),...,\chi_{\ell}(a))\tag{\nmb:{1}} 
$$ 
is an 
algebraic isomorphism. Here $\{\chi_j\}$ are the characters of the 
fundamental representations. Finally 
$$G = 
\bigcup_{a\in \on{Reg}_{\#}(G)}{\overline {\conj_G(a)}}\tag{\nmb:{2}} 
$$ is a disjoint union. 
\endproclaim 
 
\proclaim{\nmb.{4.5}. Theorem} 
Let $G$ be a reductive complex algebraic group and let 
$\pi:G\to\Aut(V)$ be a 
locally faithful rational representation of $G$. 
Let $a\in G$ be regular. Assume that $a$ is 
nonsingular with respect to the Cayley map $\Phi = \Phi_{\pi}$ so 
that $\Phi(a)$ is regular in $\g$ by \nmb!{2.5.2}. Then 
$\Phi$ restricts to an isomorphism 
$$ 
\Phi:\overline {\conj_G(a)}\to \overline {\Ad_G(\Ph(a))}\tag{\nmb:{1}} 
$$ 
of affine varieties. 
\endproclaim 
 
\demo{Proof} 
Let $G^s$ be the simply-connected covering group of the 
commutator (semisimple) subgroup $G'$ of $G$. Let 
$\ga:G^s\to G'$ be the covering map. We may write 
$a = bc$ where $b\in \on{Cent}(G)$ and $c\in G'$. Let $g\in G^s$ be such 
that $\ga(g) = c$; note that $g$ is regular. 
Clearly the mapping $G^s\to G$, given by $h\mapsto b\ga(h)$, 
restricts to a surjective $G^s$-equivariant 
morphism $\be:\conj_{G^s}(g)\to \conj_G(a)$; thus 
by continuity $\be$ also restricts
to a dominant morphism 
$$ 
\be:\overline {\conj_{G^s}(g)}\to \overline {\conj_G(a)}.\tag{\nmb:{2}} 
$$ 
But also by continuity one has a dominant morphism 
$$ 
\Phi:\overline {\conj_G(a)}\to \overline {\Ad_G(\Ph(a))}\tag{\nmb:{3}} 
$$ 
and hence if $\al$ is the composite of 
\thetag{\nmb|{2}} and \thetag{\nmb|{3}} it defines a dominant morphism 
$$ 
\al:\overline {\conj_{G^s}(g)}\to \overline {\Ad_G(\Ph(a))}\tag{\nmb:{4}} 
$$ 
But then the cohomomorphism of \thetag{\nmb|{4}} is injective. 
However since the affine 
algebras in question are $\Ad$-complete by Theorems \nmb!{4.3} and 
\nmb!{4.4}, it follows 
that \thetag{\nmb|{4}} must be an isomorphism. 
But then obviously \thetag{\nmb|{2}} and \thetag{\nmb|{3}} are isomorphisms. 
\qed\enddemo 
 
\subhead\nmb.{4.6}  \endsubhead 
Richardson proved that for semisimple 
groups the generalized Cayley map defines an isomorphism of the 
unipotent variety in G with the nilcone in $\g$. His theorem is 
very general and includes the case of $G$ defined over fields 
of finite characteristic as long as the prime is good. An 
application of Theorem \nmb!{4.5} yields Richardson's theorem for the 
complex case. 
 
\proclaim{Theorem} 
Let $a\in G$ be 
a principal unipotent element. Then $a$ is non-singular. Let 
$U\subset G$ be the unipotent variety 
$U= \overline {\conj_G(a)}$. Then 
$$ 
\Phi:U\to \overline {\Ad_G(\Ph(a))}\tag{\nmb:{1}} 
$$ 
is an isomorphism of 
affine varieties. Furthermore if $\Phi(e) = 0$ 
(e.g. $G$ is semisimple) 
then $\Phi(a)$ is principal nilpotent so that $\overline 
{\Ad_G(\Ph(a))}$ is the nilcone $N\subset \g$ and hence \thetag{\nmb|{1}} 
is an isomorphism 
$$ 
\Phi:U\to N\tag{\nmb:{2}} 
$$ 
\endproclaim 
 
\demo{Proof} 
Since $e\in U$ and 
$d\Ph(e)$ is invertible it is immediate that $d\Ph(a)$ is invertible. 
Thus \thetag{\nmb|{1}} follows from Theorem \nmb!{4.5}. 
If $\Phi(e) = 0$ then $0\in  \overline {\Ad_G(\Ph(a))}$. 
But this implies that $\Phi(a)$ is principal 
nilpotent and hence one has \thetag{\nmb|{2}}. 
\qed\enddemo 
 
\proclaim{\nmb.{4.7}. Corollary} If $G$ is semisimple and $X$ is a 
principal nilpotent element in $\g$ and if $\pi$ is irreducible, then  
$\Ph\i(X)\subset G$ consists of $|Z(G)|$ many elements. 
\endproclaim 
 
\demo{Proof} 
Let $b\in \Ph\i(X)$. Since $d\Ph(b)$ is invertible we have 
$b\in G^b=G^X=Z(G)\x (G^X)_0$. Also by \nmb!{4.6} there is a unique 
$a\in \Ph\i(X)\cap U$. 
Since $\pi$ is irreducible and $Z(G)$ 
is a finite group, for $c\in Z(G)$ we have $\pi(c)=k_c.1_V$ where  
$k_c\in S^1\subset \Bbb C$, thus 
$\tr(\pi'\Ph(c.g)\pi'(Y))=\tr(\pi(c)\pi(g)\pi'(Y)) 
=k_c.\tr(\pi(g)\pi'(Y)) 
=\tr(\pi'\Ph(g)\pi'(k_c.Y))$ 
for all $Y\in \g$. This implies $\Ph(c.g)=k_c.\Ph(g)$ for each 
$g\in G$. 
Choose the unique elements $a_c\in \Ph\i(\frac1{k_c}X)\cap U$ for all 
$c\in Z(G)$. Then $c.a_c\in G^X$ is in the coset $c.(G^X)_0$ and 
$\Ph(c.a_c)=X$. So $\Ph\i(X)$ consists of $|Z(G)|$ many elements.  
\qed\enddemo 
 
\subhead\nmb.{4.8} \endsubhead 
For any $a\in G$ (resp. $x\in \g$) let $G^a$ (resp. $G^x$) 
denote the centralizer of $a$ (resp. $x$) in $G$. Let $\g^a = 
\on{Lie}\,G^a$ (resp. $\g^x = \on{Lie}\,G^x$). If $a$ is semisimple and 
$G^a_o$ is the identity component of $G^a$ then $G^a_o$ is a 
reductive subgroup of $G$ (see Proposition 13.19 in \cit!{2} on p\. 321). 
If $x$ is semisimple then $G^x$ is connected (see e.g\. Theorem 
22.3, p\. 140 in \cit!{6}) and reductive (since $\g^x$ is clearly 
reductive). For any semisimple element $b\in G$ let $U_b$ be the 
unipotent variety in $G^b_o$ and $N_b$ the nilpotent cone in 
$\g^b$. 
Moreover, let $l = \on{rank}\,G$ and recall that an element $a\in G$ 
(resp. $x\in \g$) is called regular if $\dim\,G^a = l$ or equivalently 
$\dim\,\g^a =l$ (resp. $\dim\, G^x = l$ or equivalently 
$\dim\,\g^x = l$). 
 
For $a\in G$ with Jordan decomposition $a=a_s\,a_u$ 
let $\frak c = cent\,\g^{a_s}$ and let 
$\frak s = [\g^{a_s},\g^{a_s}]$ so that $\frak s$ is semisimple and 
$$ 
\g^{a_s}=\frak c\oplus \frak s. 
\tag{\nmb:{1}}$$

\proclaim{Proposition} 
Let $G$ be a connected reductive algebraic group. 
\roster 
\item"(\nmb:{2})" 
     Let $a\in G$ so that $a = a_s\,a_u$. Then $a_u$ is in the 
     unipotent variety $U_{a_s}$ of the reductive subgroup 
     $G^{a_s}_o$ of $G$. Conversely if $b\in U_{a_s}$ then $g = a_sb$ 
     is the Jordan decomposition of $g$ so that 
$$ 
a_sU_{a_s} =\{ b\in G\mid b_s = a_s\} 
$$ 
\item"(\nmb:{3})" 
     If $a\in G$ then $a$ is regular if and only if $a_u$ is 
     principal unipotent in  $G^{a_s}_o$. 
\item"(\nmb:{4})" 
     If $X\in \g$ then $X$ is regular if and only if $X_n$ is 
     principal nilpotent in $\g^{X_s}$. 
\item"(\nmb:{5})" 
     If $a\in G$ let $H$ be a Cartan subgroup of $G$ which contains 
     $a_s$ and let $\h = \on{Lie}\,H$. Let $C$ be the center of 
     $G^{a_s}_o$.  Then $G^{a_s}_o = C\,S$ and $a_s\in C$ so that 
     $a_s$ and hence $a$ are in $G^{a_s}_o$. Furthermore $C\subset H$ 
     and $\frak c\subset \h$ and in addition $\on{Lie} C= \frak c$. 
\endroster 
\endproclaim 
 
In spite of \thetag{\nmb|{5}} it is not true in general that $C$ is 
connected. In fact if $a_s$ corresponds to a vertex of the 
fundamental simplex then $\frak c = 0$ and $C$ is finite. 
 
\demo{Proof} \therosteritem{\nmb|{2}} 
It clearly suffices to show that 
$a_u$ is in the identity component of 
$G^{a_s}$. But this is immediate from the bijection between $U$ and 
$N$ defined by the exponential map. The latter implies that if 
$x=\log(a_u)$ then $x\in \g^{a_{s}}$. Hence the one parameter 
subgroup defined by $x$ is in $G^{a_s}_o$ so that $a_u\in 
G^{a_s}_o$. 
 
\therosteritem{\nmb|{3}} and \therosteritem{\nmb|{4}} 
It suffices by an identical 
argument to consider only the group case. Clearly 
by the uniqueness of the Jordan decomposition and 
\therosteritem{\nmb|{2}} we have 
$$ 
\g^a=\g^{a_s}\cap \g^{a_u} = (\g^{a_s})^{a_u} 
\tag{\nmb:{6}}$$ 
Of course $\on{rank}\,\g^{a_s}= l$ since 
$a_s$ is conjugate to an element in $H$. 
Thus 
$$ 
\on{rank}\,\frak c + \on{rank}\,\frak s = l=\on{rank}\,G 
\tag{\nmb:{7}}$$ 
But if $S$ is the semisimple subgroup corresponding to $\frak s$ then 
clearly  
$U_{a_s}\subset S$ so that $a_u\in S$ and hence by \thetag{\nmb|{6}}  one has 
$$ 
\g^a = \frak c \oplus \frak s^{a_u} 
\tag{\nmb:{8}}$$ 
so that $\dim\,\g^a = \on{rank}\,\frak c + \dim\,\frak s^{a_u}$. 
But $\dim\,\frak s^{a_u} \geq \on{rank}\,\frak s$ and by definition of 
principal unipotent one has $\dim\,\frak s^{a_u} = \on{rank}\,\frak s$ 
if and only if $a_u$ is principal unipotent in $S$ or equivalently in 
$G^{a_s}_o$. Now the result follows from \thetag{\nmb|{7}}. 
 
\therosteritem{\nmb|{5}} 
That $G^{a_s}_o = C\,S$ follows from \thetag{\nmb|{1}}.  
Obviously $H\subset G^{a_s}_o$ and hence $H$ is a Cartan 
subgroup of $G^{a_s}_o$. But since 
$H\subset G^{a_s}_o$ it follows that $a_s\in G^{a_s}_o$. But then 
obviously $a_s\in C$. But since the center of a connected 
reductive group lies in every Cartan subgroup one has $C\subset H$. 
We get $\on{Lie} C= \frak c$ since the centers correspond to each other 
under the Lie subgroup -- Lie subalgebra correspondence. 
\qed\enddemo 
 
\subhead\nmb.{4.9} \endsubhead 
Let $z_a = \log(a_u)$ so that $z_a\in N_{a_s}$. By the theorem of 
Jacobson-Morosov, see \cit!{7},~p\.983,  
there exists $h_a\in \frak s$ so that $[h_a,z_a] = 2z_a$. 
But then if $r_a(t) =\exp(t\,h_a)$ one has clearly has that 
$$ 
\lim_{t\to -\infty}r_a(t)\,a_u \,r_a(t)^{-1}= 1 
$$ 
But of course 
$r_a(t)$ commutes with $a_s$ since $a_s\in C$. Thus 
$$ 
\lim_{t\to -\infty}r_a(t)\,a\,r_a(t)^{-1}  = a_s 
\tag{\nmb:{1}}$$ 
 
\proclaim{\nmb.{4.10}. Corollary} 
Let $G$ be a connected reductive complex algebraic group 
and let $\Ph$ be the Cayley mapping of a rational 
representation. 
 
If for some $a\in G$ the differential 
$d\Ph(a_s):T_{a_s}G\to \g$ is invertible then also 
$d\Ph(a):T_{a}G\to \g$. This is the case if $a_s$ is hyperbolic, 
by \nmb!{5.4} below. 
\endproclaim 
 
\demo{Proof} This is an immediate consequence of \nmb!{4.9.1}. 
\qed\enddemo 
 
\subhead\nmb.{4.11} \endsubhead 
We will begin to establish results 
leading to the main theorem on the commutativity of the generalized 
Cayley mapping and the operation of taking the semisimple part for Jordan 
decompositions. We will use the notation of \nmb!{4.8}. 
 
\proclaim{Theorem} 
Let $G$ be a connected reductive complex algebraic group 
and let $\Ph$ be the Cayley mapping of a rational 
representation. Let $a\in G$. 
 
Then for the semisimple parts we have 
$\Ph(a_s)=\Ph(a)_s$ and the Jordan decomposition 
is the decomposition into components with respect to 
\nmb!{4.8.8} 
$$ 
\Ph(a) = \Ph(a_s) +\Ph(a)_n \in \g^a= \frak c \oplus \frak s^{a_u}. 
$$ 
\endproclaim

\demo{Proof} 
Let $\Ph(a)=Z+F\in \g^a= \frak c \oplus \frak s^{a_u}$ be the 
decomposition into components with respect to \nmb!{4.8.8}. 
Recall from \nmb!{4.9} the curve $r_a:\Bbb R\to  S$ 
satisfying 
$\lim_{t\to -\infty}r_a(t)\,a\,r_a(t)^{-1}  = a_s$ 
by \nmb!{4.9.1}. 
Hence by the continuity of $\Phi$ one has 
$$\align 
\Phi(a_s)&=\lim_{t\to -\infty}\Ad_{r_a(t)}(\Phi(a)) 
=\lim_{t\to -\infty}\Ad_{r_a(t)}(Z+F) 
=Z+\lim_{t\to -\infty}\Ad_{r_a(t)}F. 
\endalign$$ 
By \nmb!{2.4.2} and \nmb!{4.8.1} we have 
$\Ph(a_s)\in \on{Cent}(\g^{a_s})=\frak c$. But also 
$\Ad_{r_a(t)}F\in \frak s$ by \nmb!{4.9} so that 
$\lim_{t\to -\infty}\Ad_{r_a(t)}F=0$ which implies that $F$ is 
nilpotent, and $Z=\Ph(a_s)$ which is semisimple since 
$\frak c\subseteq \h$ by \nmb!{4.8.5}. 
Finally note that $[Z,F]=0$ since $\frak c=\on{Cent}(g^{a_s})$ 
so that the result follows. 
\qed\enddemo 
 
\subhead\nmb.{4.12} \endsubhead 
We now consider the nilpotent part of 
$\Phi(a)$. The situation is more complicated. Let 
$b\in G$ be semisimple. If 
$w\in U_b$ then of course 
$$ 
(bw)_s = b\quad\text{ and }\quad (bw)_u= w. 
\tag{\nmb:{1}}$$ 
For any $w \in U_b$ one has $\Phi(bw)_n\in N_b$ since 
$\Phi(bw)_n$ is, by the uniqueness of the Jordan 
decomposition, clearly invariant under $Ad\,b$. One thus obtains a 
map 
$$\Phi_b:U_b\to N_b,\quad \Phi_b(w) = \Phi(bw)_n. 
\tag{\nmb:{2}}$$ 
 
\proclaim{Proposition} 
Let $b\in G$ be semisimple. The map $\Phi_b:U_b \to N_b$ is a regular 
morphism which commutes with the adjoint actions of $G^b$. 
\endproclaim 
 
\demo{Proof} The commutation of $\Phi_b$ with the adjoint action of 
$G^b$ is obvious. The rationality is an immediate consequence of 
Theorem \nmb!{4.11} since it clearly imply that 
$$ 
\Phi_b(w) = \Phi(bw)-\Phi(b). \qed\tag{\nmb:{3}} 
$$ 
\enddemo 
 
\proclaim{\nmb.{4.13}. Theorem} 
Let $G$ be a connected reductive complex algebraic group 
and let $\Ph$ be the Cayley mapping of a rational 
representation with $\Ph(e)=0\in \g$. 
Let $b\in G$ be semisimple and suppose that $d\Phi(b):T_b G\to \g$ is 
invertible. This holds for $b$ hyperbolic, see \nmb!{5.4}. 
 
Then 
$\Phi_b:U_b\to N_b$ is an isomorphism of algebraic varieties. 
\endproclaim 
 
This is generalization of theorem \nmb!{4.6}. It also follows 
directly from the `lemme fondamental' in \cit!{14}, as in 
Richardsons proof of \nmb!{4.6}. 
 
\demo{Proof} The proof of \nmb!{4.6} that $\Ph:U\to N$ is an 
algebraic isomorphism is a consequence of the fact that $\Ph$ 
carries a principal unipotent orbit to a principal nilpotent 
orbit. Replacing $G$ by the semisimple algebraic group $S$ (using 
the notation of \nmb!{4.8} where 
$b=a_s$) the same argument yields the isomorphism of $\Ph_b$ as 
soon as we demonstrate that $\Ph_b$ carries a principal unipotent 
orbit in $U_b$ to a principal nilpotent orbit in $N_b$. Let $w\in 
U_b$ be a principal unipotent element in $S$. But $bw$ is a 
regular element in $G$ by \nmb!{4.8.3}. Using the notation of 
the proof of Proposition \nmb!{4.8} where $a=bw,\, a_s=b,\,a_u =w$ one has 
$\dim\,\g^{bw} = l$ and $\g^{bw} = \frak c \oplus \frak s^w$ where if $c 
= \dim\,\frak c$ and $s=\frak s^w$ then $l = c + s $ 
and $s=\on{rank}\,\frak s$. 
One the other hand if $\frak v = \frak s^{\Ph_b(w)}$ and 
$v = \dim\,\frak v$ then 
$v\geq s $ and $v=s$ if and only if $\Ph_b(w)$ is principal 
nilpotent in $S$. But, by Theorem \nmb!{4.11} clearly 
$$ 
\frak c \oplus \frak v\subset \g^{\Ph(bw)} 
\tag{\nmb:{1}}$$ 
But now by assumption $d\Ph(b)$ is invertible. 
Thus $d\Ph(bw)$ is invertible by Proposition \nmb!{4.10}. 
But then $\g^{bw} = \g^{\Ph(bw)}$ by Theorem 1.8. In particular 
$\dim\,\g^{\Ph(bw)} = l$. But the left side of \thetag{\nmb|{1}} has 
dimension $c + v$. Thus one must have $v=s$ (and equality in 
\thetag{\nmb|{1}}. Hence $\Ph_b(w)$ is a principal nilpotent in $S$. 
\qed\enddemo 
 
\head\totoc\nmb0{5}. The mapping degree of the Cayley mapping \endhead 
 
\subhead\nmb.{5.1}. The mapping degree of $\Ph$ for compact $G$ 
\endsubhead 
Let now $G$ be a compact group and the ground field be the reals. 
Then $\Ph(G)$ is compact in $\g$, so $\Ph:G\to \g$ is not surjective. 
Let us embed $\g$ into the one-point-compactification 
$\g\cup\infty = S^n$, then the topological mapping degree of 
$\Ph:G\to \g\cup\infty$ is 0 since $\Ph$ is not surjective. Thus over 
a regular value of $\Ph$ which exists by the theorem of Sard, there 
is an even number of sheets: on one half of these $\Ph$ is 
orientation preserving, on the other half it is orientation 
reversing. 
 
\subhead\nmb.{5.2}. The mapping degree in the complex case 
\endsubhead 
Since in general there are conjugacy orbits on $G$ which map to points in 
$\g$, the mapping $\Ph$ is not proper in the sense of (Hausdorff) 
topology (which means that compact sets have compact inverse images; 
in the usual topology on $\g$). But by \nmb!{2.4.1} the 
mapping $\Ph$ induces a mapping between the algebraic orbit spaces 
$\bar \Ph:  G//\conj_G\to \g//\Ad_G$, i.e\. the affine varieties with 
coordinate rings $A(G)^G$ and $A(\g)^G$, respectively. 
 
\proclaim{Theorem} 
Let $G$ be a connected reductive complex Lie group 
and let $\Ph$ be the Cayley mapping of a 
representation. 
Let $\h\subset \g$ be a Cartan subalgebra with Cartan 
subgroup $H$. 
If the Cayley mapping $\Ph:H\to \h$ is proper 
then we have:
\roster 
\item"(\nmb:{1})" 
       The mapping degree of $\Ph:H\to \h$ is positive 
       and consequently the mapping is surjective. 
\item"(\nmb:{2})" 
       $\bar \Ph:G//\conj_G\to \g//\Ad_G$ is a proper mapping for 
       the Hausdorff topologies on the affine varieties. Thus 
       $\bar\Ph$ has positive mapping degree and is 
       surjective. 
\endroster 
\endproclaim 
 
\demo{Proof} 
Let $\h$ be a Cartan subalgebra with Cartan subgroup 
$H$. By \nmb!{2.4.5} we have $\Ph:H\to\h$.
 
\therosteritem{\nmb|{1}} 
$\Ph:H\to \h$ is a proper smooth mapping; thus its mapping 
degree is defined as the value of the induced mapping in the top De Rham 
cohomology with compact supports which is isomorphic to $\Bbb R$ 
via integration: 
$$ 
\on{deg}(\Ph)=\Ph^*:H^l_c(\h;\Bbb R)=\Bbb R\to 
     H^l_c(H,\Bbb R)=\Bbb R,\quad l=\dim(\h). 
$$ 
Choose a regular value $Y$ of $\Ph$. Then for each 
$g\in \Ph\i(Y)$ the tangent mapping $T_g\Ph$ is invertible 
and orientation preserving, since it is complex holomorphic. Thus 
the mapping degree is the number of sheets over a regular point, 
which is positive. But then $\Ph:H\to\h$ has to be surjective: if 
not, its image is closed, and a $n$-form with support in the 
complement is pulled backed to 0 on $H$, in contradiction to the 
positivity of the mapping degree. 
 
\therosteritem{\nmb|{2}} 
Since $G//\conj_G\cong H/W$ and $\g//\Ad_G\cong \h/W$ where $W$ is the 
Weyl group, the result follows directly from \therosteritem{\nmb|{1}}. 
\qed\enddemo 
 
\proclaim{\nmb.{5.3}. Lemma} 
Suppose that for reals $r_i$ we have $\sum_{i=1}^N r_i=0$. 
Then 
$$ 
\sum_{i=1}^N r_ie^{r_i} \ge \frac 1{2N} \sum_{i=1}^N r_i^2. 
$$ 
\endproclaim 
 
\demo{Proof} 
We separate negative and positive summands and consider for $s_i>0$ 
and $t_j>0$ the expression 
$$ 
-\sum_{i=1}^n s_ie^{-s_i} + \sum_{j=1}^m t_je^{t_j}, 
\qquad 
\sum_{i=1}^n s_i = A = \sum_{j=1}^m t_j. 
$$ 
The function $f(x)=x\,e^x$ is strictly convex for $x>0$, thus 
$$ 
\frac1m\sum_{j=1}^m f(t_j) \ge f\Bigl(\frac1m\sum_{j=1}^m t_j\Bigr) 
=f(A/m) 
$$ 
with equality only if all $t_j$ are equal. 
Moreover, $-\sum_{i=1}^n s_ie^{-s_i}\ge -\sum_{i=1}^n s_i=-A$. For 
the right hand side we have 
$A^2=\Bigl(\sum_{j=1}^m t_j\Bigr)^2\ge \sum_{j=1}^m t_j^2$ and 
similarly $A^2\ge \sum_{i=1}^n s_i^2$. 
So finally 
$$\align 
-\sum_{i=1}^n s_ie^{-s_i} + \sum_{j=1}^m t_je^{t_j} &\ge 
-A + m\frac Am e^{A/m} 
\ge -A + A\bigl(1+\frac Am\bigr) = \frac{A^2}m = \frac1{2m}2A^2\\ 
&\ge \frac1{2m} \Bigl(\sum_{i=1}^n s_i^2 +  \sum_{j=1}^m 
t_j^2\Bigr).\qed 
\endalign$$ 
\enddemo 
 
\proclaim{\nmb.{5.4}. Theorem} 
Let $G$ be a connected reductive complex algebraic group 
and let $\Ph$ be the Cayley mapping of a rational 
representation $\pi$. 
Let $G_{\text{hyp}}\subset G$ be the subset of all semisimple 
hyperbolic elements in $G$, and let $\g_{\text{hyp}}$ be the set 
of all semisimple hyperbolic elements in the Lie algebra $\g$. 
 
Then for each $g\in G_{\text{hyp}}$ the differential 
$d\Ph(g):T_gG\to \g$ is invertible. 
 
Moreover, the mapping $\Ph_{\pi}:G_{\text{hyp}}\to \g_{\text{hyp}}$ is 
bijective and a diffeomorphism between real subvarieties if 
$\Ph_{\pi}(e)=0$, or also if the representation $\pi$ is 
self-contragredient. 
\endproclaim 
 
\demo{Proof} 
Let $\h\subset \g$ be a Cartan subalgebra with Cartan subgroup $H$. 
Let $\h=\h_{\Bbb R}+ i\h_{\Bbb R}$. 
Let $g\in G$ be a hyperbolic element which we assume to lie in 
$H$. Then $g = \exp X$ where $X\in \h_{\Bbb R}$, so 
$g\in H_{\Bbb R}$. We have $\Ph(g)\in \h_{\Bbb R}$ since 
$$ 
\tr(\pi'(\Ph(g))\pi'(Y)) = \tr(\pi(g)\pi'(Y)) \in \Bbb R 
$$ 
for all $Y\in \h_{\Bbb R}$. 
We choose a maximal compact subgroup $G_u\subset G$ such that 
$i\h_{\Bbb R}\subset \on{Lie}(G_u)=:\g_u$. Then there exists a Hermitian 
inner product on $V$ (a Hilbert space structure) such that all 
elements of $\pi(G_u)$ are unitary. Then $\pi'(\g_u)$ consists of 
skew Hermitian operators, and $\pi'(i\g_u)$ consists of Hermitian 
operators. Thus $\pi'(\h_{\Bbb R})$ and $\pi(g)$ are Hermitian 
operators. 
We use $\tr(AB)$ as real positive inner product on the the 
space of Hermitian operators on $V$. 
 
We have $g=h^2$ for $h=\exp(\tfrac12\log(g))\in H_{\Bbb R}$. 
We claim that $\pi(g):V\to V$ is a positive definite Hermitian 
operator, i.e\. $(X,Y)\mapsto \tr(\pi(g)\pi'(X)\pi'(Y)^*)$ is a 
positive definite Hermitian form on $\g$. 
Namely, for $X\in \g$ arbitrary we have $X=X_1+X_2$ for unique 
$X_1\in \g_u$ and $X_2\in i\g_u$. Let $X^*:= X_1-X_2$, then 
$\pi'(X^*)$ equals the adjoint operator $\pi'(X)^*$. We have 
$$\multline 
\tr(\pi(g)\pi'(X)\pi'(X)^*) = 
\tr(\pi(h)^2\pi'(X)\pi(X)^*) =\\ 
=\tr(\pi(h)\pi'(X)(\pi(h)\pi'(X))^*) =\|\pi(h)\pi'(X)\|^2>0. 
\endmultline$$ 
By \nmb!{2.4.7} this implies that 
$d\Ph(g):T_gG\to \g$ is invertible. 
 
For the second assertion, 
we claim that in both cases we have for $C>0$ 
$$ 
\|\pi'(\Ph(g))\| \ge C\|\pi'(X)\|. 
\tag{\nmb:{1}}$$ 
 
If $\Ph(e)=0$ we have 
$$\align 
0 &= \tr(\pi'(\Ph(e))\pi'(X)) = \tr(\pi(e)\pi'(X)) = \tr(\pi'(X)) \\ 
&= \sum_{\la\in\text{weight}(\pi)} \tr(\pi'(X)|V_\la) 
     = \sum_{\la\in\text{weight}(\pi)} \dim(V_\la)\la(X). \\ 
\endalign$$ 
Thus by lemma \nmb!{5.3} we have 
$$ 
\sum_{\la\in\text{weight}(\pi)} \dim(V_\la)\la(X)e^{\la(X)} \ge 
C\sum_{\la\in\text{weight}(\pi)} \dim(V_\la)\la(X)^2 
$$ 
for a positive constant $C$. 
Then we have by Cauchy Schwarz, 
$$\align 
\|\pi'\Ph(g)\|\|\pi'(X)\| &\ge \tr(\pi(g)\pi'(X)) 
     = \sum_{\la\in\text{weight}(\pi)} 
     \tr(\pi(\exp(X))\pi'(X)|V_\la) \\ 
&= \sum_{\la\in\text{weight}(\pi)} 
     \dim(V_\la)\la(X)e^{\la(X)} 
     \ge C\sum_{\la} \la(X)^2\dim(V_\la) \\ 
&= C\sum_{\la} \tr(\pi'(X)^2|V_\la) 
     =C\tr(\pi'(X)^2) = C\|\pi'(X)\|^2. 
\endalign$$ 
If the representation $\pi$ is 
self-contragredient, a similar argument works using 
$re^r-re^{-r}\ge 2r$ for $r>0$ instead of the inequality \nmb!{5.3}. 
 
Now we may finish the proof. 
Note that $\exp:\h_{\Bbb R}\to H_{\Bbb R}$ is a diffeomorphism with 
inverse $\log:H_{\Bbb R}\to \h_{\Bbb R}$. Thus estimate \thetag{\nmb|{1}} 
implies that $\Ph:H_{\Bbb R}\to\h_{\Bbb R}$ is a proper mapping. 
It is also 
a local diffeomorphism, thus a covering mapping and a diffeomorphism since 
$\h_{\Bbb R}$ is vector space. 
 
Finally, each hyperbolic element $g\in G$ is contained in 
$H_{\Bbb R}$ for a suitable Cartan subgroup $H$, and the above 
arguments show that $\Ph:G_{\text{hyp}}\to \g_{\text{hyp}}$ is 
locally a diffeomorphism and is surjective. It is also injective: Let 
$\Ph(g_1)=Y=\Ph(g_2)$. Then $g_1$ and some conjugate 
$h\,g_2\,h\i$ lie in the same Cartan subgroup $H_{\Bbb R}$ on which 
$\Ph$ is a diffeomorphism, thus $g_1=h.g_2.h\i$. By equivariancy 
we have $\Ad(h)Y=Y$. Since $\Ph$ is a local diffeomorphism near 
$g_1$, the orbits have the same dimension, thus $g_1$ and $Y$ have 
the same connected component of the isotropy group. But isotropy 
groups of semisimple elements of the Lie algebra are connected, as 
mentioned in \nmb!{4.8}.
Thus $g_1=g_2$. 
\qed\enddemo 
 
\proclaim{\nmb.{5.5}. Theorem} 
Let $G$ be a connected reductive complex algebraic group 
and let $\Ph$ be the Cayley mapping of a rational 
representation with $\Ph(e)=0\in \g$. 
 
Then $\Ph:G_{\text{pos}}\to \g_{\text{real}}$ is bijective and a 
fiber respecting isomorphism of real algebraic varieties, where 
$G_{\text{pos}}$ and $\g_{\text{real}}$ are defined in \nmb!{4.1}. 
\endproclaim 
 
\demo{Proof} 
Let $g\in G_{\text{pos}}$, then $g_s=g_h$, thus by \nmb!{5.4} 
$d\Ph(g_s):T_{g_s}G\to \g$ is invertible and thus by theorem 
\nmb!{4.13} $\Ph_{g_s}:U_{g_s}\to N_{g_s}$ is an automorphism of 
algebraic varieties, where $U_{g_s}$ is the unipotent variety in the 
reductive group $G^{g_s}_o$, see \nmb!{4.8}. 
The set 
$G_{\text{pos}}= 
\bigsqcup_{h\in G_{\text{hyp}}}h.U_h\to G_{\text{hyp}}$ is a 
fibration with complex algebraic varieties as fibers and the real 
algebraic variety $G_{\text{hyp}}$ as base. Likewise 
$\g_{\text{real}}= 
\bigsqcup_{h\in G_{\text{hyp}}}(\Ph(h)+N_h)\to \g_{\text{hyp}}$ is a 
is a fibration with the nilpotent cones $N_h$ as fibers. 
$\Ph:G_{\text{hyp}}\to \g_{\text{hyp}}$ is given by 
$\Ph(g)=\Ph(g_s)+\Ph_{g_s}(g_u)=\Ph(g)_s+\Ph(g)_n$ and is a fiber 
respecting isomorphism by theorem \nmb!{4.13} and \nmb!{5.4}. 
\qed\enddemo 
 
\head\totoc\nmb0{6}. Examples \endhead 
 
\subhead\nmb.{6.1}. The Cayley mappings for the representations of 
$SL_2(\Bbb C)$ 
\endsubhead 
 
The standard representation $SL_2(\Bbb C)\subset \End(\Bbb C^2)$. Here 
$\Ph(A)=A-\frac12\tr(A)1_{\Bbb C^2}$ and $\Ps(A)=\frac12\tr(A)$ for 
$A\in SL_2$. $\Ph$ is surjective and proper and has mapping degree 2. 
On the Cartan subgroup we get 
$$ 
\Ph_{\pi_1}\pmatrix a & 0 \\ 0 & a\i \endpmatrix = 
\pmatrix \tfrac{1-a^2}{2a} & 0 \\ 0 & \tfrac{a^2-1}{2a} \endpmatrix, 
\qquad \Ps_{\pi_1}\pmatrix a & 0 \\ 0 & a\i \endpmatrix = \frac{a^2+1}{2a}. 
$$ 
 
For the $(n+1)$-dimensional representation we computed
$$\align 
\Ph_{\pi_n}\left(\smallmatrix a & 0 \\ 0 & a\i \endsmallmatrix\right) 
&= \frac{3}{n^3+3n^2+2n} 
     \sum_{p=0}^n (n-2p) a^{n-2p} 
\left(\smallmatrix 1 & 0 \\ 0 & -1 \endsmallmatrix\right) \\
\Ps_{\pi_n}\left(\smallmatrix a & 0 \\ 0 & a\i \endsmallmatrix\right) 
&=\frac{-108}{(n^3+3n^2+2n)^3} 
\sum_{p_1,p_2,p_3=0}^n (n-2p_1)^2(n-p_2)n  
   a^{3n-2(p_1+p_2+p_3)} 
\endalign$$ 
 
\proclaim{\nmb.{6.2}. Example: The standard representation of 
$SL_n(\Bbb C)$} 
The standard representation $SL_n(\Bbb C)\subset \End(\Bbb C^n)$. Here 
$\Ph(A)=A-\frac1n\tr(A)1_{\Bbb C^n}$, and $\Ph$ has mapping degree 
$n$, which coincides with the degree of the smooth hypersurface 
$SL_n\subset \End(\Bbb C^n)$. 
Thus $\Ph$ maps regular orbits in $SL_n$ to 
regular orbits in $\frak s\frak l_n$, and by \nmb!{3.2} we have 
$$ 
A(SL_n) = A(SL_n)^{SL_n}\otimes\on{Harm}(SL_n)). 
$$ 
 
The set $G_{\text{sing}}$ where the differential of the Cayley 
transform is singular is given by 
$\Ps\i(0)=(SL_n\cap \frak s\frak l_n)\i = \{A\in SL_n:\tr(A\i)=0\}$. 
 
Let $\ch\in A(SL_n)$ be given by $\ch(A)=\frac1n\tr(A)$. Then 
$$\align 
A(SL_n)^{SL_n}&=A(\frak s\frak l_n)^{SL_n}[\ch],\\ 
Q(SL_n)^{SL_n}&=Q(\frak s\frak l_n)^{SL_n}[\ch],\\ 
Q(SL_n)&=Q(\frak s\frak l_n)[\ch], 
\endalign$$ 
in terms of \nmb!{3.3}. 
\endproclaim 
 
This can be shown as follows: 
$\End(\Bbb C^n)=\frak s\frak l_n\oplus \Bbb C.1_n$ is an orthogonal 
direct sum, thus the orthogonal projection is given by 
$\pr(A)=A-\frac1n\tr(A)$. Of course $\Ph=\pr|SL_n$. 
For $Y\in\frak s\frak l_n$ we have $A\in \Ph\i(Y)$ if and only if 
$A-Y=t.1_n$ for some $t\in \Bbb C$. The set of all these $t$ is given 
by the equation $\det(Y+t.1_n)=1$. Generically there are $n$ 
solution, thus the degree of $\Ph$ is $n$, and $\Ph$ is surjective. 
 
For $A\in SL_n$ and $X\in \frak s\frak l_n$ then $A.X$ is a typical 
tangent vector in $T_A(SL_n)$.  
Now $d\Ph(A)=\Ph:T_A(SL_n)\to \frak s\frak l_n$ is not invertible if 
its kernel $\Bbb C.1_n\cap T_A(SL_n)$ is non trivial, so there exists 
an $X\in \frak s\frak l_n$ with $A.X=1_n$ (by scaling $X$ 
appropriately). But then $X=A\i\in \frak s\frak l_n$ and thus 
$\tr(A\i)=0$. 
 
Let us discuss the adjunction now. The fiber over any 
$X\in \frak s\frak l_n$ consists of all $A=X+t.1_n$ where $t$ runs 
through the roots of the polynomial $f_X(t)=\tr(t.1_n+X)-1$, and 
moreover $t=\ch(A)$, by the formula for $\Ph$. Let us consider the 
expansion 
$$ 
f_X(t) = \tr(t.1_n+X)-1 = \sum_{j=0}^np_j(X)\, t^j, 
$$ 
where $p_j\in A(\frak s\frak l_n)^{SL_n}$ form a system of 
generators with $p_0(X)=\det(X)-1$, $p_{n-1}(X)=\tr(X)=0$, and 
$p_n(X)=1$. We have 
$$ 
f_X(\ch) = \tr(t.1_n+X)-1 = \sum_{j=0}^np_j(X)\, \ch^j = 0, 
$$ 
and this is the minimal polynomial for $\ch$ over 
$A(\frak s\frak l_n)^{SL_n}$. Finally note that $\ch$ is $1/n$ times 
the character of the standard representation and $\Ph^*$ pulls back 
the elements $p_j\in A(\frak s\frak l_n)^{SL_n}$ to the characters of 
the remaining fundamental representations of $SL_n$. Thus 
$A(SL_n)^{SL_n}=A(\frak s\frak l_n)^{SL_n}[\ch]$. The other 
assertions follow from \nmb!{3.3}. 
 
\subhead\nmb.{6.3}. The standard representation of $O_n(\Bbb C)$ 
\endsubhead 
In the standard representation $O_n(\Bbb C)\to \End(\Bbb C^n)$ we 
have 
$$\align 
O_n(\Bbb C)&=\{A\in \End(\Bbb C^n): AA^\top = 1_n\}\\ 
\frak s\frak o_n(\Bbb C)&=\{X\in \End(\Bbb C^n):X+X^\top=0\}\\ 
\Ph(A) &= \tfrac12(A-A^\top). 
\endalign$$ 
 
\head\totoc\nmb0{7}. Spin representations and Cayley transforms \endhead 
 
In this section we treat only the complex group 
$\on{Spin}(n,\Bbb C)$. Some of the results here were inspired by 
\cit!{13}. 
We first recall notations and results from \cit!{12}. 
 
\subhead\nmb.{7.1}. Clifford multiplication in terms of exterior 
algebra operations 
\endsubhead 
Let $V$ be a complex vector space of finite dimension $n>0$. 
Assume that $B_V=(\quad,\quad)$ is a fixed nonsingular symmetric bilinear 
form on $V$. Its extends naturally to a symmetric nonsingular 
bilinear form $B_{\wedge V}$ on the exterior algebra $\wedge V$. 
The graded structure on $\wedge V$ induces a $\Bbb Z$-graded structure 
on the operator algebra $\End(\wedge V)$ and also on the 
$\Bbb Z$-graded super Lie algebra $\on{Der}(\wedge V)$ of graded 
super derivations. 
For $\beta\in \End(\wedge V)$ we denote its transpose by 
$\beta^t\in \End(\wedge V)$ so that $(\beta u,v)= (u,\beta^tv)$ for 
any $u,v\in\wedge V$. 
 
For any $u\in \wedge V$ let $\epsilon(u)\in \End(\wedge V)$ (left 
exterior multiplication) be defined so that $\epsilon(u)w =u\wedge w$ 
for any $w$. For any $u\in \wedge V$ let $\iota(u)= \epsilon(u)^t$. 
If $u\in \wedge^kV$ then $\epsilon(u)\in \End^k(\wedge V)$ and 
$\iota(u)\in  \End^{-k}(\wedge V)$. Furthermore if $x\in V$ then 
$\iota(x)\in \on{Der}^{-1} \wedge V$. In fact $\iota(x)$ is the 
unique element of $\on{Der}^{-1}(\wedge V)$ such that $\iota(x)y = 
(x,y)$ where we have identified $\Bbb C = \wedge^0V$. For a basis 
$\{z_j\}$ of $V$ we have $\on{Der}(\wedge V) = 
\bigoplus_j\epsilon(\wedge V)\iota(z_j)$. Let 
$\kappa\in \End^0(\wedge V)$ be defined so that $\kappa = (-1)^k$ on 
$\wedge^k V$. Let $\wedge^{\text{even}}V$ and $\wedge^{\text{odd}}V$ 
be the eigenspaces for $\kappa$ corresponding, respectively, to the 
eigenvalues $1$ and $-1$. 
 
We recall that the Clifford algebra $C(V)$ over $V$ with respect to 
$B_V\vert V$ is the tensor algebra $T(V)$ over $V$ modulo the ideal 
generated by all elements in $T(V)$ of the form $x\otimes x -(x,x)$ 
where $x\in V$ and we regard $\Bbb C = T^0(V)$. We note that 
$\epsilon(x)\iota(y) + \iota(y)\epsilon(x) = (x,y) 
\in \End(\wedge V)$. For any $x\in V$ let $\gamma(x) = \epsilon(x) + 
\iota(x)$ then $\gamma(x)^2 = (x,x)$. Thus the correspondence 
$x\mapsto \gamma(x)$ extends to a homomorphism 
$C(V)\to \End(\wedge V)$, denoted by $u\mapsto \gamma(u)$, defining 
the structure of a $C(V)$-module on $\wedge V$. This leads to the map 
$u\mapsto \gamma(u) 1$, denoted by  $\psi:C(V)\to \wedge V$.  By 
Theorem II.1.6, p. 41 in \cit!{3} this map $\psi$ is bijective and we 
will identify $C(V)$ with $\wedge V$, using $\psi$. We consequently 
then recognize that the linear space $\wedge V = C(V)$ has 2 
multiplicative structures. If $u,v\in \wedge V$ there is exterior 
multiplication $u\wedge v$ and Clifford multiplication $uv = 
\gamma(u)v$. Also we will write $C^0(V) = \wedge^{\text{even}}V$ and  
$C^1(V) = \wedge^{\text{odd}}V$ so that 
$C^i(V)C^j(V)\subset C^{i+j}(V)$ for $i,j\in \Bbb Z/(2)$.  That is, 
the parity automorphism $\kappa$ of the exterior algebra is also an 
automorphism of the Clifford algebra. In addition if 
$\alpha\in \End \wedge V$ is the unique involutory antiautomorphism 
of $\wedge V$, as an exterior algebra, such that $\alpha(x) = x$ for 
all $x\in V$, then $\alpha$ is also an antiautomorphism of the 
Clifford algebra structure; explicitly, $\alpha = 
(-1)^{{k(k-1)\over 2}}$ on $\wedge^kV$. Using $\al$ and $\ka$ one can 
show easily that for $x\in V$ and $u\in C(V)$ we have 
$ux=(\ep(x)-\io(x))\ka(u)$. 
 
 
Let $\on{Spin}(V)$ be the set of all $g\in C^0(V)$ such that 
$g\alpha(g) =1$ and $gx\alpha(g)\in V$ for all $x\in V$. This is an 
algebraic group under Clifford multiplication. For any 
$g\in \on{Spin} V$ let $T(g)\in GL(V)$ be the mapping given by 
$T(g)x = gx\alpha(g)$ for $x\in V$. 
 
The Lie algebra $\on{Lie}\on{Spin} V$ is a Lie subalgebra of 
$C^0(V)$, in fact $\on{Lie}\on{Spin} V= \wedge^2V$ as explained on 
p.~68 in \cit!{3}. We need this in detail: Let $u\in \wedge^2V$. For 
any $x\in V$ one has $\iota(x)u\in V$ and one can define 
$\tau(u)\in \frak s\frak o(V)\subset\End(V)$ by $\tau(u)x = 
-2\iota(x)u$. 
 
For any $z\in \End V$ there exists a unique operator 
$\xi(z)\in \on{Der}^0(\wedge V)$ which extends the action of $z$ on 
$V$. Clearly $\xi:\End V\to  \on{Der}^0(\wedge V)$ is a Lie algebra 
isomorphism. 
 
\proclaim{\nmb.{7.2}. Theorem} 
$\on{Spin} V\subset C^0(V)$ is a connected Lie group. If 
$g\in \on{Spin} V$ then $T(g)\in SO(V)$ and 
$T:\on{Spin} V\to SO(V)$ is an epimorphism with kernel $\{\pm 1\}$. 
In particular $T$ defines $\on{Spin} V$ as a double cover of 
$SO(V)$. 
 
The subspace $\wedge^2V$ is a Lie subalgebra of $C^0(V)$. In fact 
$\wedge^2V = \on{Lie}\on{Spin} V$ so that $\on{Spin} V$ is the group 
generated by all $\exp(u)$ for $u\in \wedge^2V$ where exponentiation 
is with respect to Clifford multiplication. 
 
Furthermore the map $\tau:\wedge^2V\to \on{Lie}SO(V)$ is a Lie 
algebra isomorphism and is the differential of the double cover 
$T:\on{Spin} V\to SO(V)$. Finally if $u\in \wedge^2V$ then $\ad(u) = 
\xi(\tau(u))$ so that not only is $\ad(u)$ a derivation of both 
exterior algebra and Clifford algebra structures on $C(V)$ but also 
$\ad:\wedge^2V\to \on{Der}^0\wedge V$ is a faithful representation of the Lie 
algebra $\wedge^2V$ on $\wedge V$. 
\endproclaim 
 
\subhead\nmb.{7.3} \endsubhead 
Let $\{z_i\}$ be a $B_V$-orthonormal basis of $V$. 
On these elements the Clifford product equals the exterior product, 
$z_iz_j=z_i\wedge z_j$. 
Let $N =\{1,\ldots,n\}$ and let ${\Cal I}$ be the set of all subsets 
of $N$. We regard any subset $I$ as ordered: $I = \{i_1<\ldots<i_k\}$. 
Let 
$$z_I = z_{i_1}\cdots z_{i_k} = z_{i_1}\wedge \cdots\wedge  z_{i_k}.$$ 
The set of elements $\{z_I\}$ with $\vert I\vert  = k$, is a 
$B_{\wedge V}$-orthonormal basis of $\wedge^kV$ so 
that the  set $\{z_I\}, I\in {\Cal I}$, $\vert I\vert$ is even (resp. odd) 
is a 
$B_{\wedge V}$-orthonormal basis of $C^0(V)$ (resp. $C^1(V)$) and the full set 
$\{z_I\},\,I\in {\Cal I}$ is a $B_{\wedge V}$-orthonormal basis of $C(V)$. 
Moreover for $I,J\in {\Cal I}$ let 
$I\diamond J\in {\Cal I}$ be the symmetric difference 
$(I\setminus J)\cup (J\setminus I)$. Then 
$$\gather 
z_Iz_J = c_{I,J}z_{I\diamond J},\qquad 
\text{ where }c_{I,J}\in\{1,-1\},\\ 
(z_I)^{-1} = \alpha(z_I) = (-1)^{\vert I\vert (\vert I\vert-1)/2} z_I. 
\endgather$$ 
 
\subhead\nmb.{7.4}. The spin module \endsubhead 
The following description of the spin module $S$ is uniform in 
$n$ and the spin representation will always be a direct sum of 2 
(possibly equivalent) irreducible representations. 
 
Let $Z = \wedge^0V + \wedge^nV$ so that $Z$ is a 2-dimensional 
subspace of $C(V)$ and let $C^Z(V)$ be the centralizer of $Z$ in 
$C(V)$. There exists (uniquely up to sign) an element 
$\mu\in \wedge^nV$ such that $\mu^2 =1$.  Hence if we put $e_+ = 
{1\over 2}(1 +\mu)$ and  $e_- = {1\over 2}(1 -\mu)$ then 
$\{e_+,e_-\}\subset Z$ are orthogonal idempotents in the sense of ring 
theory and $1 = e_+ + e_-$. In particular, $Z$ is a subalgebra of 
$C(V)$ and is isomorphic to $e_+Z\oplus e_-Z= 
\Bbb C\oplus \Bbb C\subset C^Z(V)$. We have $x\mu = (-1)^{n+1}\mu v$ 
for any $v\in V$ so that $C^0(V)\subset C^Z(V)$ and in fact $C^Z(V) = 
ZC^0(V)$. Let $C^Z_+(V) = e_+\,C^Z(V)$ and $C^Z_-(V) = e_-\,C^Z(V)$. 
Let $n_1 = [{n+1\over 2}]$. It is well known (see e.g\. II.2.4 and 
II.2.6 in \cit!{3}) that both $C_+^Z(V)$ and $C_-^Z(V)$ are each 
isomorphic to a $2^{n_1-1}\times 2^{n_1-1}$ matrix algebra so that 
$C^Z(V)$ is a semisimple associative algebra of dimension 
$2^{2n_1-1}$ and the unique decomposition of $C^Z(V)$ into simple 
ideals is given by $C^Z(V) = C_+^Z(V) \oplus C_-^Z(V)$. In particular 
$Z= \on{Cent}(C^Z(V))$ and $e_+$ and $e_-$ are, respectively, the identity 
elements of $C_+^Z(V)$ and $C_-^Z(V)$. Thus there exists a 
$C^Z(V)$-module $S$ (the spin module) $S$ of dimension $2^{n_1}$, 
defined by a representation $ \sigma:C^Z(V)\to \End S$ characterized 
uniquely, up to equivalence, by the condition that if $S_+ = e_+\,S$ 
and $S_- = e_-\,S$ then $S= S_+\oplus S_-$ is the unique 
decomposition of $S$ into irreducible $C^Z(V)$-modules.  Let 
$\sigma_+:C^Z(V)\to \End S_+$ and $\sigma_-:C^Z(V)\to \End S_-$ be the 
corresponding irreducible representations. 
 
Recall $\on{Spin}(V)\subset C^0(V)\subset C^Z(Z)$. 
The restriction of $\sigma$ 
to $\on{Spin}(V)$ will be denoted by $\on{Spin}$ so that $\on{Spin}$ 
is a representation of $\on{Spin}(V)$ on $S$. Its differential, also 
denoted by $\on{Spin}$, is the restriction $\sigma$ to $\wedge^2V$. 
Thus one has a Lie algebra representation 
$\on{Spin}:\wedge^2V\to \End S$. Replacing $\sigma$ by $\sigma_{\pm}$ 
and $S$ by $S_{\pm}$ one similarly has representations 
$\on{Spin}_{\pm}$ of $\on{Spin} V$ and $\wedge^2V$ on $S_{\pm}$. 
Since $\wedge^2V$ generates $C^0(V)$, both $\on{Spin}_{\pm}$ are 
irreducible representations. 
 
Finally, let  $n_0 = [{n\over 2}]$ so that one has $n = n_0 + n_1$. 
Note that if $n$ is even then $n_1 =n_0 = {1\over 2}n$, whereas if 
$n$ is odd then $n_1 = n_0 +1$ so that $n_1= {1\over 2}(n+1)$ and 
$n_0 = {1\over 2}(n-1)$. 
 
\subhead\nmb.{7.5} \endsubhead 
Let $\pr_k:C(V)\to \wedge^k V$ be the projection defined by 
the graded structure of $\wedge V$. We identify $\wedge^0V$ with $\Bbb C$ 
so that the image of $\pr_0$ is scalar-valued. 
 
The Clifford algebra $C(V)$ as a module over itself by left 
multiplication $(\gamma)$ is equivalent to $2^{n_0}$ copies of the 
$C(V)$ spin module $S$ with respect to $\sigma$ (see \nmb!{7.4}). One 
has that $2^n = 2^{n_0}\,2^{n_1}$ 
 
\proclaim{Proposition} 
For any $w\in C(V)$ one has 
$$ 
\pr_0(w)= {1\over 2^n}  \tr  \gamma(w) 
={1\over 2^{n_1}} \tr  \sigma(w). 
$$ 
\endproclaim 
 
\demo{Proof} 
Of course for $c\in \wedge^0 V$ one obviously has that $c= {1\over 
2^n} \tr  \gamma(c)$. With the notation of \nmb!{7.3}, 
it suffices to prove that $\tr \gamma(z_I)=0$ if $k>0$. 
But this follows from \nmb!{7.3}. 
\qed\enddemo 
 
\subhead\nmb.{7.6} \endsubhead 
We now note  that the bilinear form given by $\tr \sigma$ 
is essentially given by $B_{\wedge V}$. 
 
\proclaim{Proposition} 
Let $u,w\in C(V)$ Then 
$$ 
\tr \sigma(u)\sigma(w) = 2^{n_1}(u,\alpha(w)) 
$$ 
\endproclaim 
 
\demo{Proof} The left side 
is just $\tr \sigma(uw)$. By proposition \nmb!{7.5} this is just 
$2^{n_1}\pr_0(uw)$. Using the basis $z_I$ of \nmb!{7.3} 
and writing $u = \sum_Iu_Iz_I$ and 
$w = \sum_Iw_Iz_I$ it is immediate from \nmb!{7.3} that 
$\pr_0(uw) = \sum_Iu_I w_I \pr_0(z_I^2)$. 
But clearly by \nmb!{7.3} one has $(z_I,\alpha(z_J)=0$ for $I\neq J$ and 
$z_I^2=\pr_0(z_I^2)= (z_I,\alpha(z_I)$ . 
\qed\enddemo 
 
\subhead\nmb.{7.7}. The Cayley mapping \endsubhead 
The generalized Cayley mapping 
$\Ph_T:\on{Spin}(V)\to \on{Lie}\on{Spin}(V)$ for the representation 
$T:\on{Spin}(V)\to SO(V)$ is given by 
$\tr(\sigma(g)\sigma(y)) = 
\tr(\sigma(\Ph_T(g))\sigma(y))$. By \nmb!{7.2} one has 
$\on{Lie}\on{Spin}(V) = \wedge^2V$. 
 
\proclaim{Proposition} 
Let $\g\in \on{Spin}(V)$. Then $\Ph_T(g) = \pr_2(g)$ 
\endproclaim 
 
\demo{Proof} Let $y\in \wedge^2V$. Let 
$x= \Ph_T(g)$ and let $z=\pr_2(g)$. 
Note that $\alpha(y) = - y$. By \nmb!{7.6} one has 
$$ 
\tr(\sigma(g)\sigma(y)) = -2^{n_1}(g,y)= -2^{n_1}(\pr_2(g),y) 
= \tr(\sigma(\pr_2(g))\sigma(y)).\qed 
$$ 
\enddemo 
 
\subhead\nmb.{7.8}\endsubhead 
Let $\theta:SO(V)\to \on{Aut}(\wedge V)$ be the representation 
so that if $a\in SO(V)$ then $\theta(a)$ is the unique exterior algebra 
automorphism which extends the action of $a$ on $V$. Clearly 
$$\tr \theta(a) = \det (1+a)$$ 
 
\proclaim{Proposition} 
For any $g\in \on{Spin}(V)$ and $w\in C(V)$ then 
$$ 
\theta(T(g)) w = g\,w g^{-1} 
$$ 
In particular for $\theta(a)$, for $a\in SO(V)$ is 
an automorphism of both the Clifford and exterior algebra 
structures on $C(V)$. 
\endproclaim 
 
\demo{Proof} 
Using the notation of \nmb!{7.1} and \nmb!{7.2} 
it is clear that $\xi\vert \on{Lie}SO(V)$ is the differential of 
$\theta$. To prove the proposition it suffices to establish its 
infinitesimal analogue. But this is stated in \nmb!{7.2}. 
\qed\enddemo 
 
\proclaim{\nmb.{7.9}. Proposition} 
The tensor product representation 
$\on{Spin}\otimes \on{Spin}$ of $\on{Spin}(V)$ on $S\otimes 
S$ descends to $SO(V)$ and is equivalent to $\theta$ if $n$ is even 
and 2 copies of $\theta$ if $n$ is odd. 
\endproclaim 

This is well known; we include a proof in the conventions used above. 

\demo{Proof} Assume that $M$ is a matrix algebra and $\beta$ is an 
antiautomorphism of $M$. Let $L$ be a minimal left ideal of $M$ so 
that $(\dim L)^2 = \dim M$. Then clearly $R  = \beta(L)$ is a minimal 
right ideal of $M$. 
There exists $v\in M$ such that $LvR\neq 0$ so that $Lv R = M$ 
since $LvR$ is a 
2-sided ideal and $M$ is simple. This implies that the map 
$$ 
\mu:L\otimes L \to M,\qquad \mu(a\otimes b)= a\,v\,\beta(b) 
$$ 
is surjective and thus a linear isomorphism by 
dimension. If $g\in M$ then clearly 
$\mu(ga\otimes gb)= g\mu(a\otimes b) \beta(g)$, which by proposition 
\nmb!{7.8} proves the assertion  
in case $n$ is even by choosing $M= C(V)$, 
$\beta = \alpha$, $L=S$ and $g\in \on{Spin}(V)$. 
 
If $n$ is odd then 
$C(V) = M_1\oplus M_2$ where $M_i,\,i=1,2,$ are matrix algebras. Then 
$S = L_1 \oplus L_2$ where $L_i,\,i=1,2,$ are respectively minimal 
left ideals in $M_i$ (and left ideals in $C(V)$). Then 
$$ 
S\otimes S = L_1\otimes L_1 \oplus L_1\otimes L_2 
\oplus L_2\otimes L_1 \oplus L_2\otimes L_2 
\tag{\nmb:{1}}$$ 
We have 
$\on{Spin}(V)\subset C^0(V)$ which implies that all 4 summands on 
the right side of \thetag{\nmb|{1}} define $\on{Spin}(V)$-equivalent 
submodules of 
$\on{Spin}\otimes \on{Spin}:\on{Spin}(V)\to \on{Aut}(S\otimes S)$, 
since the parity automorphism $\kappa\in \on{Aut} C(V)$ from \nmb!{7.1} 
of course fixes $C^0(V)$ but interchanges $M_1$ and $M_2$. 
This is clear from 
the definition of $e_{\pm}$ in \nmb!{7.4} since $\kappa = -1$ on 
$\wedge^nV$. Now if $i\in \{1,2\}$ let $i'\in \{1,2\}$ be defined so 
that $\alpha(M_i) = M_{i'}$. But then $\alpha(L_i)=R_{i'}$ is a 
minimal right ideal (and right ideal in $C(V)$) in $M_{i'}$. But then 
$$ 
C(V) = L_1\, M_1\,R_{1'}\oplus L_1\,M_1\,R_{2'}\ 
\oplus L_2\,M_2\,R_{1'}\oplus L_2\,M_2\,R_{2'}, 
$$ 
where  
2 of the 4 components are identically zero and the remaining 
2 are $M_1$ and $M_2$. Note that the 2-sided ideal $M_i$ in $C(V)$ is 
stable under $\theta$ by proposition \nmb!{7.8}. The argument in the even 
case then readily implies that $\theta$ is equivalent to 2 of the 4 
components in \thetag{\nmb|{1}} and hence $\on{Spin}\otimes \on{Spin}$ is 
equivalent to 2 copies of $\theta$. \qed\enddemo 
 
\proclaim{\nmb.{7.10}. Theorem} 
Let $g\in \on{Spin}(V)$ then $\pr_0(g) =0$ if 
$g\notin \on{Spin}(V)^*:= \{g\in \on{Spin}(V) : \det (1+T(g))\neq 0\}$. 
On the 
other hand if $g\in \on{Spin}(V)^*$ then 
$$ 
\pr_0(g) = {1\over 2^{n/2}} 
\sqrt {\det (1+ T(g))} 
$$ 
for one of the two choices of 
the square root. Both choices are taken for the two elements 
$g,g'\in \on{Spin}(V)$ such that 
$T(g)=T(g')$. 
\endproclaim 
 
\demo{Proof} 
Apply Proposition \nmb!{7.5} to the 
case where $w=g\in \on{Spin}(V)$. Then 
$$ 
\pr_0(g)^2 = {1\over 2^{2n_1}}(\tr \on{Spin}(g))^2. 
$$ 
But $(\tr \on{Spin}(g))^2= 
\tr (\on{Spin}\otimes \on{Spin})(g)$. If $n$ is even then $2n_1 = n$ and 
we have 
$\tr (\on{Spin}\otimes \on{Spin})(g) = \det (1+T(g))$ by \nmb!{7.8} 
and proposition \nmb!{7.9}. 
Thus in this case 
$$ 
\pr_0(g)^2 = {1\over 2^{n}} \det (1+T(g)). 
\tag{\nmb:{1}}$$ 
Now assume 
that $n$ is odd. Then $\tr (\on{Spin}\otimes \on{Spin})(g) = 2  
\det (1+T(g))$ by \nmb!{7.8} and Proposition \nmb!{7.9}. 
But since $2n_1-1 = n$ equation \thetag{\nmb|{1}} holds 
for all $n$. The theorem now follows from: 
If $g_1,g_2\in \on{Spin}(V)$ then 
$T(g_1) = T(g_2)$ if and only if $g_1=\pm g_2$. 
\qed\enddemo 
 
\subhead\nmb.{7.11} \endsubhead 
For any $u\in \wedge^2V$ let $\hat{e}^u$ denote the 
exponential of $u$ using exterior (and hence 
commutative - since $\wedge^{\text{even}}V$ is commutative) multiplication. 
If 
$w\in 
\wedge^2V$ and $x\in V$ then by \nmb!{7.1} and \nmb!{7.2} one has, by 
differentiating an exponential, 
$$ 
-\iota(x)\hat{e}^{2w} = \varepsilon([w,x]){e}^{2w}, 
$$ 
where the bracket is Clifford commutation. But 
then also 
$$ 
(-\varepsilon([w,x])-\iota(x))\hat{e}^{2w} = 
(\iota(x)+\varepsilon([w,x]))\hat{e}^{2w}. 
$$ 
Adding 
$(\varepsilon(x) + \iota([w,x]))\hat{e}^{2w}$ to both sides of the 
last equation 
yields 
$$ 
(\varepsilon(x-[w,x])- \iota(x-[w,x]))\hat{e}^{2w} = (\varepsilon(x+[w,x])+ 
\iota(x+[w,x]))\hat{e}^{2w}. 
$$ 
But this is just Clifford 
multiplication. If $y\in V$ and $u\in C(V)$ then 
$(\varepsilon(y) +\iota(y))u = yu$ and, if $u\in 
C^0(V)$, 
$(\varepsilon(y) -\iota(y))u = uy$, both by \nmb!{7.1}. 
We have proved: 
 
\proclaim{Proposition} 
Let $x\in V$ and $w\in \wedge^2V$. Then 
$$ 
\hat{e}^{2w}(x-[w,x]) = (x+ [w,x])\hat{e}^{2w}. 
$$ 
\endproclaim 
 
\subhead\nmb.{7.12} \endsubhead 
Let $s\in \on{Lie}SO(V)$ 
and let $V_s^1$ and $V_s^{-1}$ be, respectively, the $s$-eigenspaces 
for the eigenvalues $\pm 1$. For $x,y\in V$ the equation $$(sx,y) = 
-(x,sy)$$ readily implies 
 
\proclaim{Lemma} 
The subspaces $V_s^1$ and 
$V_s^{-1}$ are isotropic with respect to the symmetric bilinear form 
$B_V$ and are nonsingularly paired by $B_V$. In particular 
$$ 
\dim V_s^1= \dim V_s^{-1} 
$$ 
Consequently $1+s$ 
is invertible if and only if $1-s$ is invertible. 
\endproclaim 
 
\subhead\nmb.{7.13} \endsubhead 
Using lemma \nmb!{7.12} we define the Zariski open subset of 
$\on{Lie}SO(V)$ and $SO(V)$ by 
$$\align 
\on{Lie}SO(V)^{(*)}&=\{s\in \on{Lie}SO(V):\dim(V^1_s)=\dim(V^{-1}_s)=0\}\\ 
&=\{s\in \on{Lie}SO(V):1+s\text{ and }1-s\text{ are invertible}\},\\ 
SO(V)^* &= \{a\in SO(V)\mid \det (1+a)\neq 0\}, 
\endalign$$ 
so that in 
the notation of theorem \nmb!{7.2} and theorem \nmb!{7.10} one has 
$$ 
\on{Spin}(V)^* = T^{-1}(SO(V)^*). 
$$ 
 
\proclaim{Proposition} 
The subsets $\on{Lie}SO(V)^{(*)}$ and $SO(V)^*$ are algebraically 
isomorphic via the mapping
$$ 
\on{Lie}SO(V)^{(*)}\to SO(V)^*,\qquad s\mapsto  a= {1-s\over 1+s}, 
\qquad s = {1-a\over 1+a} \gets a
\tag{\nmb:{1}}$$ 
for $s\in \on{Lie}SO(V)^{(*)}$ and $a\in SO(V)^*$.
In addition one has the relation $(1+a)(1+s) = 2$. 
\endproclaim 
 
\demo{Proof} 
If $b\in \End V$ let $b^t$ be the transpose endomorphism with respect to 
$B_V$. Since $s^t= -s$ for $s\in \on{Lie}SO(V)$ it is immediate from 
\thetag{\nmb|{1}} 
that $a^t = a^{-1}$ so that $a$ is orthogonal. Since 
$\on{Lie}SO(V)^{(*)}$ is Zariski open it is connected and $0\in 
\on{Lie}SO(V)^{(*)}$. It follows then that $a\in SO(V)$. Furthermore 
$(1+a)(1+s) = 2$ is immediate from \thetag{\nmb|{1}} and hence $a\in 
SO(V)^*$. The map $s\mapsto a$ is injective since one recovers $s$ from 
$a$ by \thetag{\nmb|{1}}. Conversely if $a\in SO(V)^*$, 
then $s\in \on{Lie}SO(V)$ since $a^t = a^{-1}$ 
implies $s^t = -s$. But this formula yields $(1+a)(1+s) = 2$ so that $s\in 
\on{Lie}SO(V)^{(*)}$. Hence the mapping \thetag{\nmb|{1}} 
is bijective. 
\qed\enddemo 
 
\subhead\nmb.{7.14} \endsubhead 
In general if $b\in \End V$ is such that $1+b$ is invertible we will 
put 
$$ 
\Ga(b) = {1-b\over 1+b} 
$$ Also let 
$\wedge^2V^{(*)}$ be the inverse image of $\on{Lie}SO(V)^{(*)}$ 
under the isomorphism $\tau:\wedge^2V\to \on{Lie}SO(V)$. 
 
\proclaim{Theorem} 
Let $g\in \on{Spin}(V)^*$. Let the sign of square root be chosen so that 
$$ 
\pr_0(g) = {1\over 2^{n/2}} \sqrt {\det (1+ T(g))} 
$$ 
(see 
Theorem \nmb!{7.10}). Then putting $c = \pr_0(g)$ one has 
$$ 
g = c\,\,\hat {e}^{-2w} 
$$ 
where $w =\ta\i(\Ga(T(g)))\in \wedge^2V$. 
In particular $w\in \wedge^2V^{(*)}$. 
\endproclaim 
 
\demo{Proof} We have 
$\Ga(\tau(w)) = T(g)$. 
But by proposition \nmb!{7.11}, for any $x\in V$, one has 
$\hat{e}^{-2w}((1+\tau(w))x) = ((1-\tau(w))x)\hat{e}^{-2w}$. 
Putting $y = (1+\tau(w))x)$ one therefore has 
$$ 
\hat{e}^{-2w}\,\, y = T(g)(y)\,\,\hat{e}^{-2w} 
$$ 
On the other hand 
$g\,y\,g^{-1} = T(g)\,(y)$ so that 
$$ 
y\,g^{-1} = g^{-1}T(g)(y) 
$$ 
But the last two equations imply that 
$g^{-1}\,{e}^{-2w}$ commutes with $y$. Since $y\in V$ is arbitrary 
this implies that $g^{-1}\,{e}^{-2w} \in \on{Cent}(C(V))\cap C^0(V)$. 
But $\on{Cent}(C(V))\cap C^0(V) = \Bbb C$. Hence there exists $d\in 
\Bbb C$ such that $d\, {e}^{-2w} = g$. But by applying $\pr_0$ it 
follows that $d=c$. 
\qed\enddemo 
 
\subhead\nmb.{7.15} \endsubhead 
We determine the generalized Cayley mapping 
$\Ph:\on{Spin}(V)\to \on{Lie}\on{Spin}(V)$, corresponding to to the 
representation $\on{Spin}$, for elements $g$ in the Zariski 
open set $\on{Spin}(V)^*$. 
 
\proclaim{Theorem} 
Let $g\in \on{Spin}(V)^*$. Let the sign of square root be chosen so that 
$$\pr_0(g) = {1\over 2^{n/2}}\,\sqrt {\det (1+ T(g))},$$ 
according to theorem \nmb!{7.10}. Then 
$$\Ph_{\on{Spin}}(g)= -2\,\pr_0(g)\,\ta\i(\Ga(T(g)))\in \wedge^2V.$$ 
Thus the generalized Cayley mapping 
$\Ph_{\on{Spin}}:\on{Spin}(V)\to \wedge^2V$  
factors to the classical 
Cayley transform $\Ga:SO(V)^*\to \on{Lie}\on{Spin}(V)^{(*)}$, up to
multiplication by a regular function,
via the natural identifications. 
 
For the degree of the spin representation we have 
$$\deg(\on{Spin}) =\cases n \quad &\text{ if }n \text{ is even,}\\ 
                n-1\quad&\text{ if }n \text{ is odd.}\endcases 
$$ 
Let $\ch\in A(\on{Spin}(n,\Bbb C))$ be given by 
$\ch(a)=\sqrt{\det(1+T(a))}=2^{n/2}\pr_0(a)$. Then 
$$\align 
Q(\on{Spin}(n))^{\on{Spin}(n)}&=Q(\on{Lie}\on{Spin}(n))^{\on{Spin}(n)}[\ch],\\ 
Q(\on{Spin}(n))&=Q(\on{Lie}\on{Spin}(n))[\ch], 
\endalign$$ 
in terms of \nmb!{3.3}. 
\endproclaim 
 
\demo{Proof} 
By \nmb!{7.7} one has $\Ph_{\on{Spin}}(g) = \pr_2(g)$ 
but since the exponential in  theorem \nmb!{7.14} is an exponential 
for exterior multiplication the result follows from this theorem.

Let us now determine $\deg(\on{Spin})$. Consider the  fiber of $\Ph$ over a 
generic $X$ in $\on{Lie}\on{Spin}(n)$. Let $G_o = 
\{a\in \on{Spin}(n)\mid \det(1+T(a)) = 0\}$. We may assume 
that $X$ is in the complement of $\Phi(G_o)$. Then modulo some fixed 
scalar (which can  be ignored 
in determining the cardinality of the fiber) if $a\in\on{Spin}(n)$ is 
in the fiber then 
$$ 
\sqrt{\det (1+T(a))}(1-T(a))(1+T(a))\i =X. 
$$ 
Put $t=\sqrt{\det(1+T(a))}$ so that $(1-T(a))(1+T(a))\i=X/t$. 
Recall that an equation of matrices $(1-c)(1+c)\i=d$ 
is symmetric in $c$ and $d$ and that $(1+c)(1+d)=2$. 
Thus 
$$ 
(1+T(a))(1+X/t)=2.\tag{\nmb:{1}} 
$$ 
Taking the determinant of both sides of \thetag{\nmb|{1}}  one has 
$$ 
t^{2}\det(1+X/t)=2^n\tag{\nmb:{2}} 
$$ 
or 
$$ 
f_X(t)=\det(t+X)-2^nt^{n-2}=0. \tag{\nmb:{3}} 
$$ 
Comparing \thetag{\nmb|{1}} and \thetag{\nmb|{2}} we get 
$\det(1+T(a))=t^2$ which checks with the definition of $t$. 
Given a non-zero root $t$ of \thetag{\nmb|{3}} the corresponding 
$T(a)$ is given by 
$T(a)=(1-X/t)(1+X/t)\i$. But $t$ is a specific square root of 
$\det(1+T(a))$ and hence $a$ itself is uniquely defined. 
 
Now \thetag{\nmb|{3}} is a polynomial of degree $n$ in $t$ with 
leading term $t^n$ and constant term $\det(X)$. If $n$ is even then 
$\det(X)\ne 0$ for generic $X$ and we also may assume that the polynomial 
$\det(t+X)$ has pairwise different roots. Replacing $X$ by a large 
multiple $rX$ and considering 
$$ 
\tfrac1{r^n}f_{rX}(t)=\tfrac1{r^n}\left(\det(t+rX)-2^nt^{n-2}\right) 
=\det(\frac tr+X)-\tfrac{2^n}{r^2}(\tfrac tr)^{n-2} 
$$ 
which approaches $\det(t/r+X)$ we see that also the polynomial 
\thetag{\nmb|{3}} has pairwise distinct roots and they are all 
non-zero since $\det(X)\ne0$. 
 
If $n$ is odd then $\det(X)=0$ but generically the rank of $X$ is 
$n-1$; by the argument above \thetag{\nmb|{3}} has again pairwise 
distinct roots for generic $X$, but one of them is $0$ so that there 
are $n-1$ different non-zero roots. 
 
Let us discuss the adjunction now. We use the above description of the 
fiber of $\Ph$ over a generic $X$. Let us consider the expansion 
$$ 
f_X(t) = \det(t+X)-2^nt^{n-2}= \sum_{j=0}^np_j(X)\, t^j, \tag{\nmb:{4}} 
$$ 
where $p_j\in A(\on{Lie}\on{Spin}(n))^{\on{Spin}(n)}$ form a system of 
generators with $p_0(X)=\det(X)$ which is 0 for odd $n$, and 
$p_n(X)=1$. We have 
$$ 
f_X(\ch) = \det(\ch+X)-2^nt^{n-2} = \sum_{j=0}^np_j(X)\, \ch^j = 0, 
$$ 
since $t$ in the discussion above corresponds to $\ch(a)$, 
and this is the minimal polynomial for $\ch$ over 
$Q(\on{Lie}\on{Spin}(n))^{\on{Spin}(n)}$ for $n$ even. For $n$ odd 
the minimal polynomial is $f_X(\ch)/\ch$. The reason for this is that 
$f_X(t)$ has exactly $\deg(\on{Spin})$ many pairwise nonzero roots for 
generic $X$, in a Zariski dense open set which we may determine as 
the complement to $\Ph(\{g:\Ps(g)=0\})$. 
\qed\enddemo 
 
\subhead\nmb.{7.16}. The Cayley transform for the 
$\rh$-representation of any semisimple Lie group 
\endsubhead 
Let $G$ be any simply connected semisimple Lie group. 
One has a homomorphism 
$$ 
G\to\on{Spin}(\g) \tag{\nmb:{1}} 
$$ 
which lifts the adjoint representation. The Lie algebra of 
$\on{Spin}(\g)$ is $\wedge^2\g$. 
Thus for any $a\in G$ there exists $u=u(a)\in\wedge^2\g$ where 
$$ 
u = \sqrt{\det (1 + Ad a)}(1 - Ad a)/(1+ Ad a) \tag{\nmb:{2}} 
$$ 
On the other hand $\wedge\g$ is a differential chain complex with 
respect to a boundary operator $c$ where $c:\wedge^2\g\to\g$ 
is given by  $c(X\wedge Z)=[X,Z]$ for $X,Z\in \g$. 
 
\proclaim{Theorem} Let $G$ be any simply connected semisimple Lie 
group. For the $\rh$-representation one has 
up to a fixed scalar multiple 
$$ 
\Phi(a)=c(u). \tag{\nmb:{3}} 
$$ 
\endproclaim

\demo{Proof} 
Actually the first reference establishing  that the restriction 
(using \thetag{\nmb|{1}}) 
of the spin representation $s$ of $\on{Spin}(\g)$ to $G$ is a 
multiple of the $\rho$-representation is in reference 9 in \cit!{12} 
(see top paragraph on p.~358 where 
$\rh$ has been written as $g$).  But our result \nmb!{2.8.1} 
on direct sums (here a multiple of the same representation, thus we 
may use \nmb!{2.8.1} for $G$ semisimple) we may use the restriction 
$\pi=s|G$ to compute $\Phi_{\rho}$. But then 
$\pi':\g \to \wedge^2\g$ is an injection. 
Let $\frak m$ be the orthocomplement of $\pi'(\g)$ in $\wedge^2\g$ 
so that one has a direct sum 
$\wedge^2\g = \pi'(\g) + \frak m$ with corresponding projection 
$p:\wedge^2\g\to\pi'(\g)$. 
Let $a\in G$ and put $X=\Phi(a)$. 
In the notation of \thetag{\nmb|{3}} we must prove that $X = c(u)$ 
up to a scalar multiple. But one clearly has 
$$ 
\pi'(X) = p(u) \tag{\nmb:{4}} 
$$ 
That is for some $v\in\frak m$ one has 
$$ 
u = v + \pi'(X) \tag{\nmb:{5}} 
$$ 
But now besides the boundary operator $c$ on $\wedge\g$ one has a coboundary 
operator $d$ on $\wedge\g$. See section 2 in \cit!{10} where 
$\partial$ in \cit!{10} is 
written as $c$ here. One has $d:\g\to\wedge^2\g$.  The homomorphism 
$\pi'$ is $\delta$ in (69) of \cit!{12}. By (106) in \cit!{12} one  has 
$\pi'=d/2$, so that \thetag{\nmb|{4}} 
becomes 
$$ 
d(X)/2  = p(u) \tag{\nmb:{6}} 
$$ 
But by (94) in \cit!{12} one has that $dc+cd$ is 1/2 times the 
Casimir operator on $\wedge\g$. But the Casimir is normalized so that 
it takes the value 1 on $\g$. 
However $c$ vanishes on $\g=\wedge^1\g$. Thus $cd=1/2$ on $\g$. Thus upon 
applying $c$ to \thetag{\nmb|{6}} yields 
$$ 
X/4 = c(pu). \tag{\nmb:{7}} 
$$ 
On the other hand $\frak m$ is the kernel of $c$ on $\wedge^2\g$ by 
(4.4.4) in \cit!{10}. 
Thus $c$ vanishes on $v$ in \thetag{\nmb|{5}} and by applying $c$ to 
both sides of \thetag{\nmb|{5}} one has 
$c(u)=c(pu)$ and hence $X/4=c(u)$. 
But this is just \thetag{\nmb|{3}} up to a scalar. 
\qed\enddemo

\enddocument